\documentclass[a4paper]{article}

\usepackage{amsmath,amssymb,graphicx}
\usepackage{subfigure}

\usepackage{lscape}
\usepackage{mathrsfs}

\usepackage{algorithmic}
\usepackage[ruled]{algorithm}

\usepackage{multirow}
\usepackage{verbatim}
\usepackage{enumerate}
\usepackage{color}

\usepackage{longtable}

\title{ On the M-eigenvalues of elasticity tensor and the strong ellipticity condition }

\author{  Hua Xiang$^a$\thanks{%Corresponding author (H. Xiang).
E-mail: hxiang@whu.edu.cn. H. Xiang is supported by the National
Natural Science Foundation of China under grants 11571265, 11471253 and NSFC-RGC %(China-Hong Kong)
No.11661161017. }
\quad Liqun Qi$^b$\thanks{
E-mail: maqilq@polyu.edu.hk.  L. Qi is supported by the Hong Kong Research Grant Council
    (Grant No. PolyU  15302114, 15300715, 15301716 and 15300717). }
\quad Yimin Wei$^{c}$\thanks{
E-mail: ymwei@fudan.edu.cn. Y. Wei is supported by the National Natural
Science Foundation of China under grant 11771099 and International Cooperation Project of Shanghai Municipal Science and Technology Commission under grant 16510711200.
}
\\ \\
\small{$^a$ School of Mathematics and Statistics, Wuhan University, Wuhan, 430072, P.R. China}\\
{\small $^b$ Department of Applied Mathematics, The Hong Kong Polytechnic University,
Hong Kong}\\
{\small $^c$ School of Mathematical Sciences and }
\\ {\small Shanghai Key Laboratory of Contemporary Applied Mathematics, }
\\ {\small Fudan University, Shanghai, 200433, P. R. of China}
}

\begin{document}

%\date{}
\maketitle

\begin{abstract}

Strong ellipticity is an important property in the elasticity theory.  In 2009, M-eigenvalues were introduced for the elasticity tensor.  It was shown that M-eigenvalues are invariant under coordinate system choices, and the strong ellipticity condition holds if and only if all the M-eigenvalues of the elasticity tensor are positive.   Thus, M-eigenvalues are some intrinsic parameters of the elasticity tensor.    In this paper, we show that the M-eigenvalues of the elasticity tensor are closely related with some elastic moduli, such as the bulk modulus, the shear modulus, Lam\'e's first parameter, the P-wave modulus, etc, and the positiveness of the M-eigenvalues are corresponding to some existing conditions for strong ellipticity in some special cases, such as the isotropic case, the cubic case, the polar anisotropic case and the tetragonal case. We also present new sufficient conditions for the strong ellipticity of the orthotropic case.  These, in a certain sense, further reveal the physical meanings of M-eigenvalues.

\end{abstract}

%\begin{keywords}
%\end{keywords}

%\begin{AMS}
%\end{AMS}

%\tableofcontents

\section{Introduction}

The elasticity tensor $\mathcal{C}$ is a fourth-rank tensor, namely a four-linear map $\mathcal{C}$: $\mathbb{R}^3 \times \mathbb{R}^3 \times \mathbb{R}^3 \times \mathbb{R}^3 \rightarrow \mathbb{R} $. With respect to the basis $\{\partial_i \}$ in a local coordinate system, the components of $\mathcal{C}$ is expressed by $c_{ijkl} = \mathcal{C}(\partial_i, \partial_j, \partial_k, \partial_l)$, where $i,j,k,l\in \{1,2,3\}$. With respect to the dual basis $\{dx^i\}$, the elasticity tensor can be written as $\mathcal{C}= c_{ijkl} dx^i dx^j dx^k dx^l$.
Using these components, we can reformulate the tensor $\mathcal{C}$ as
$$ \mathcal{C}(w,x,y,z) = \mathcal{C} wxyz = c_{ijkl}w^i x^j y^k z^l ,$$
where $w=w^i \partial_i$, $x=x^i \partial_i$, $y=y^i \partial_i$, $z=z^i \partial_i$,  and we adopt the Einstein convention of summing over the repeated indices.
When working on the orthonormal Cartesian coordinate system with a Euclidean metric given by the Kronecker delta, the covariant and contravariant types do not differ, and the indices can be upper or lower, and hence we have the option to work with only subscripts in the remainder of this article.

For the fourth-order elasticity or stiffness tensor $\mathcal{C}$,
the strong ellipticity is defined by the following function to be positive.
\begin{equation*}
f(x,y)=\mathcal{C} xyxy = C_{ijkl} x_i y_j x_k y_l > 0, \ \ \forall  x,y \in \mathbb{R}^3.
\end{equation*}
%Here we use the Einstein summation notation.
Such strong ellipticity condition ensures that the governing differential equations for elastostatics problems be completely elliptic. It is an important property in the elasticity theory associated with uniqueness, instability, wave propagation, etc, and has been studied extensively \cite{Bower_book, SChirita_JE07, GourgiotisBigoni16, Gow17, Gurtin73, HanDaiQi09,KnowlesSternberg75, KnowlesSternberg76, Mur13, PichPri16, QiChenChen18, Rosakis90, Sfy11, SimpsonSpector, Str13, WaltonWilber, WangAron96, ZubovRudev16}.

In 2009, Qi, Dai and Han \cite{QDH09} introduced the M-eigenvalues for the elasticity tensor. Here, the letter ``M'' means mechanics.
It was shown that the strong ellipticity condition holds if and only if the optimal value of a global polynomial optimization problem is positive, which further naturally leads to the M-eigenvalue problem from the optimality condition.

The M-eigenvalues $\theta$ of the elasticity tensor $\mathcal{C}$ are defined as follows.
\begin{eqnarray}
 \mathcal{C} \cdot yxy &=& \theta x , \label{eqn:MeigDefx}\\
 \mathcal{C} ~ xyx \cdot &=& \theta y, \label{eqn:MeigDefy}
\end{eqnarray}
under the constraints $x^T x=1 $ and $y^T y=1$.
Here,  $(\mathcal{C} \cdot yxy)_i = C_{ijkl}y_jx_ky_l$ and $(\mathcal{C} ~ xyx \cdot)_l = C_{ijkl}x_iy_jx_k$ in the Einstein summation notation.
It was shown that M-eigenvalues are invariant under coordinate system choices, and the strong ellipticity condition holds if and only if all the M-eigenvalue of the elasticity tensor $\mathcal{C}$ are positive \cite{QDH09}. Thus, M-eigenvalues are some intrinsic parameters of the elasticity tensor.

The M-eigenvalue problem can be reduced to the Z-eigenvalue problem.
\begin{equation}\label{eqn:ZeigDef}
 \mathcal{C} x^3 = \eta  x ,
\end{equation}
where $x^T x=1$, and $\eta$ is the Z-eigenvalue. As the M-eigenvalue is associated with the strong ellipticity, the Z-eigenvalue is related to the rank-one positive definiteness of $\mathcal{C}$. All the Z-eigenvalues  are M-eigenvalues of $\mathcal{C}$, but not vice versa. Obviously, the positiveness of Z-eigenvalues provides the necessary conditions for the strong ellipticity of material.

We note that there exists the major symmetry $C_{ijkl}=C_{klij}$, and the minor symmetry
$C_{ijkl}=C_{jikl}$ and $C_{ijkl}=C_{ijlk}$.   The requirement of the symmetry reduces the number of different elements to 21.
%The requirement of the symmetry of the stress and strain tensors lead to equality of many of the elastic constants, reducing the number of different elements to 21.
%This number can be further reduced for some special cases.
%For example, a polar anisotropic (or transversely isotropic) material exhibits hexagonal symmetry, and the number of independent constants in the elasticity tensor are reduced to 5.
%
Checking the invariance of $\mathcal{C}$ under a specific transformation $Q$ can further reduces this number for some special cases.
If the material response is to be the same, i.e., $C_{ijkl} = Q_{im} Q_{jn} Q_{kp} Q_{lq} C_{mnpq}$, then it yields the material symmetry relation, which will provide a system of equations that allows reduction in the number of independent elastic moduli \cite{Cowin_QJMAM90,PaoloVannucci_book18}.
For example, as to an orthotropic material, it has three mutually perpendicular symmetry planes, and this type of material has nine independent material constants.
Common examples of such materials include wood and fiber-reinforced composites.
Omitting the zero elements, the M-eigenvalue problem \eqref{eqn:MeigDefx}-\eqref{eqn:MeigDefy} for this case reads

{\tiny
\begin{eqnarray*}%\label{eqn:Meig4PolarAniso}
C_{1111}y_1 x_1 y_1 + (C_{1122}+C_{1221})y_1x_2y_2 + (C_{1133}+C_{1331})y_1x_3y_3 + C_{1313}y_3x_1y_3 + C_{1212}y_2x_1y_2  &=& \theta x_1 , \\
(C_{2211}+C_{2112})y_2x_1y_1  + C_{2222}y_2x_2y_2 + (C_{2233}+C_{2332})y_2x_3y_3  + C_{2323}y_3x_2y_3  + C_{2121}y_1x_2y_1  &=& \theta x_2 ,\\
(C_{3311}+C_{3113})y_3x_1y_1   + (C_{3322}+C_{3223})y_3x_2y_2 + C_{3333}y_3x_3y_3 + C_{3232}y_2x_3y_2  + C_{3131}y_1x_3y_1 &=& \theta x_3, \\
C_{1111}x_1 y_1 x_1 + (C_{2211}+C_{1221})x_2y_2x_1 + (C_{3311}+C_{1331})x_3y_3x_1 + C_{3131}x_3y_1x_3 + C_{2121}x_2y_1x_2  &=& \theta y_1 , \\
(C_{1122}+C_{2112})x_1y_1x_2  + C_{2222}x_2y_2x_2 + (C_{3322}+C_{2332})x_3y_3x_2  + C_{3232}x_3y_2x_3  + C_{1212}x_1y_2x_1  &=& \theta y_2 ,\\
(C_{1133}+C_{3113})x_1y_1x_3   + (C_{2233}+C_{3223})x_2y_2x_3 + C_{3333}x_3y_3x_3 + C_{2323}x_2y_3x_2  + C_{1313}x_1y_3x_1 &=& \theta y_3.
\end{eqnarray*}
}

To keep accordance with the standard shorthand notation in the theory of elasticity, we need the Voigt notation or the Kelvin notation. Since the minor symmetries are valid, the fourth-order stiffness tensor
$ C_{ijkl}$  may be written as a $6 \times 6$ stiffness matrix $ ( c_{\alpha \beta } ) $
 (a tensor of second order).
 %[The constitutive tensor of linear elasticity: its decompositions, Cauchy relations, null %Lagrangians, and wave propagation]
Due to the major symmetry, this matrix is also symmetric.
In the Voigt notation, each pair of indices is compressed into one index, and the standard mapping for tensor indices reads as follows.
$$ \begin{matrix}ij & = \\ \Downarrow & \\ \alpha  & = \end{matrix} \quad \begin{matrix}11&22&33&23,32&13,31&12,21\\ \Downarrow &\Downarrow &\Downarrow &\Downarrow &\Downarrow &\Downarrow &\\1&2&3&4&5&6\end{matrix} $$
For example, in a properly chosen coordinate system, the Voigt matrix of the orthotropic material has the following form.
$$ ( c_{\alpha \beta } ) = {\begin{bmatrix}c_{11} &c_{12} & c_{13}&0&0&0 \\c_{12}&c_{22}&c_{23}&0&0&0 \\ c_{13}&c_{23}&c_{33}&0&0&0 \\0&0&0&c_{44}&0&0\\0&0&0&0&c_{55}&0\\0&0&0&0&0&c_{66}\end{bmatrix}}.$$
The elasticity tensor is described by nine parameters:
$c_{11}= C_{1111}$, $c_{22} = C_{2222}$, $c_{33}=C_{3333}$, $c_{23}=C_{2233}$, $c_{13}=C_{1133}$, $c_{12}=C_{1122}$, $c_{44}=C_{2323}=C_{3232}$, $c_{55}=C_{1313}=C_{3131}$, and $c_{66}=C_{1212}=C_{2121}$.

%****************************
Considering the symmetries in the elasticity tensor $\mathcal{C}$, we naturally introduce the following eigenproblem.
\begin{equation}\label{eqn:eigentensor}
C_{ijkl} x_{kl} = \zeta ~ x_{ij} ,
\end{equation}
where the eigentensor $x_{ij}$ is symmetric.
Such eigenproblem is closely related to the positive definiteness of $\mathcal{C}$, a less general hypothesis than the strong ellipticity, which also guarantees the uniqueness of solutions in problems of elasticity \cite[P.159]{HetnarskiIgnaczak_Book2011}.
The elasticity stiffness tensor $\mathcal{C}$ must be positive definite \cite[P.48]{PaoloVannucci_book18}\cite[P.34]{Ting_Book96}. That is, The strain energy density or elastic potential satisfies
$$ \frac{1}{2} C_{ijkl} \epsilon_{ij} \epsilon_{kl} > 0, $$
where $\epsilon_{ij}$ is any symmetric strain tensor. It physically means that energy is needed to deform an elastic body from its unloaded  equilibrium position \cite{Hoger_JE95}.  The condition for the positiveness of strain energy can be replaced by the less restrictive condition $C_{ijkl} \gamma_{ij} \gamma_{kl} > 0$, where the $\gamma_{ij}$ need not be symmetric \cite[P.34]{Ting_Book96}.

Considering the symmetry of $C$, we rewrite \eqref{eqn:eigentensor} in its equivalent matrix form.
{\tiny
$$\begin{bmatrix}
C_{1111} &C_{1122} &C_{1133} &\sqrt{2} C_{1123} &\sqrt{2} C_{1113} &\sqrt{2} C_{1112}  \\
C_{2211} &C_{2222} &C_{2233} &\sqrt{2} C_{2223} &\sqrt{2} C_{2213} &\sqrt{2} C_{2212} \\
C_{3311} &C_{3322} &C_{3333} &\sqrt{2} C_{3323} &\sqrt{2} C_{3313} &\sqrt{2} C_{3312} \\
\sqrt{2} C_{2311} &\sqrt{2} C_{2322} &\sqrt{2} C_{2333} &2 C_{2323} &2 C_{2313} &2 C_{2312} \\
\sqrt{2} C_{1311} &\sqrt{2} C_{1322} &\sqrt{2} C_{1333} &2 C_{1323} &2 C_{1313} &2 C_{1312} \\
\sqrt{2} C_{1211} &\sqrt{2} C_{1222} &\sqrt{2} C_{1233} &2 C_{1223} &2 C_{1213} &2 C_{1212}
\end{bmatrix}
\begin{pmatrix}
x_{11} \\ x_{22} \\x_{33} \\ \sqrt{2}  x_{23} \\ \sqrt{2}  x_{13} \\ \sqrt{2}  x_{12}
\end{pmatrix}
= \zeta
\begin{pmatrix}
x_{11} \\ x_{22} \\x_{33} \\ \sqrt{2}  x_{23} \\ \sqrt{2}  x_{13} \\ \sqrt{2}  x_{12}
\end{pmatrix} .
$$
}
The positiveness condition corresponds to impose
that all the eigenvalues $\zeta$ are positive \cite{PaoloVannucci_book18}.
What's more, the positive definiteness of $\mathcal{C}$ implies strong ellipticity, while the converse statement is not true \cite[P.159]{HetnarskiIgnaczak_Book2011}. In other words, the positive definiteness conditions are sufficient for the strong ellipticity.

In the following, we consider five further simplified cases of the orthotropic material: the isotropic case, the cubic case, the polar anisotropic case, the tetragonal case and the orthotropic case, which have two, three, five, six and nine independent elastic parameters respectively.
%
%In the isotropic case, we show that there are two M-eigenvalues: the first M-eigenvalue is the shear modulus, and the second M-eigenvalue is the P-wave modulus.   Since the shear modulus must be positive, this implies that strong ellipticity holds in the isotropic case if and only the P-wave modulus is positive.  In the cubic case and the polar anisotropic case, we show that the positiveness condition of M-eigenvalues is equivalent to existing condition for strong ellipticity in the literature.
%
In this paper, we will show that the M-eigenvalues are closely related with some elastic moduli, such as the bulk modulus, the shear modulus, Lam\'e's first parameter, the P-wave modulus, etc. In fact, all the elastic constants in the diagonal of stiffness matrix $( c_{\alpha \beta} )$ are M-eigenvalues. The positiveness of the M-eigenvalues are corresponding to some existing conditions for strong ellipticity in some special cases, such as the isotropic case, the cubic case and the polar anisotropic case.  Thus, in a certain sense, these further reveal the physical meanings of M-eigenvalues.

\section{Isotropy}

An isotropic material has uniform properties independent on its orientation, and only two material constants are required to characterize such material.
The elasticity tensor of isotropic material is given as follows.
\begin{equation*}
C_{ijkl}=K \delta_{ij} \delta_{kl}+G \left( \delta_{ik} \delta_{jl}+\delta_{il} \delta_{jk}-\frac{2}{3}\delta_{ij}\delta_{kl} \right),
\end{equation*}
where
$\delta_{ij}$ is the Kronecker delta, $K$ is the bulk modulus (or incompressibility), and
 $G$ is the shear modulus (or rigidity), also usually denoted by  $\mu$.
Besides, Young's modulus $E$ and Poisson’s ratio $\nu$ are the most common properties used to characterize elastic solids, where $\nu = \frac{E}{2\mu}-1$.
Young's modulus $E$ is usually large.
Poisson's ratio $\nu$ typically ranges from 0.2 to 0.49 and is around 0.3 for most
metals. For a stable material, $-1<\nu<0.5$, and the solid is incompressible if $\nu=0.5$.
The bulk modulus $K$ is usually bigger than $E$, while the shear modulus is usually somewhat smaller than $E$ \cite[p.75]{Bower_book}.
Such practical considerations will be used in Section 4 to set up the conditions for strong ellipticity.

Using the Einstein summation convention of summing on repeated indices, we can rewrite  \eqref{eqn:MeigDefx} as $C_{ijkl}y_j x_k y_l = \theta x_i$ ($i=1,2,3$). Direct calculation yields that
\begin{eqnarray}\label{eqn:MeigCx}
 C_{ijkl}y_j x_k y_l  &=& K y_i x_k x_k +G(y_i x_i y_l +y_k x_k y_i- \frac{2}{3}y_i x_k y_k)  \nonumber\\
 &=& \left( K+\frac{G}{3} \right) x_k y_k y_i +G y_l y_l x_i = \theta x_i , \quad (i=1,2,3).
\end{eqnarray}
Similarly from \eqref{eqn:MeigDefy}, we have $ C_{ijkl} x_i y_j x_k = \theta y_l$ ($l=1,2,3$).  We can calculate that
\begin{eqnarray}\label{eqn:MeigCy}
C_{ijkl} x_i y_j x_k  &=& K x_i y_i x_l +G(x_k y_l x_k +x_l y_k x_k -\frac{2}{3} x_i y_i x_l ) \nonumber\\
&=& \left( K+\frac{G}{3} \right) x_i y_i x_l +G x_k x_k y_l = \theta y_l , \quad (l=1,2,3).
\end{eqnarray}

Notice that \begin{math} x_k x_k =y_l y_l=1 \end{math}, and \begin{math} x_i y_i \equiv \langle x,y \rangle \end{math}, we can reformulate \eqref{eqn:MeigCx} and \eqref{eqn:MeigCy} as
\begin{equation*}
\left( K+\frac{G}{3} \right) \langle x,y \rangle y_i=(\theta-G)x_i, \qquad  i=1,2,3. \end{equation*}
\begin{equation*}
\left( K+\frac{G}{3} \right) \langle x,y \rangle x_l=(\theta-G)y_l, \qquad l=1,2,3.
\end{equation*}

We can see that
\begin{enumerate}[(i)]
\item If $\langle x, y \rangle = 0$, then $\theta = G$.
\item If not, then $\theta \neq G$. We can derive that $x = \pm y$, and $K+\frac{G}{3}=\theta-G$, that is, $\theta = K+\frac{4}{3}G$.
\end{enumerate}
Thus, we have two M-eigenvalues:
$$\theta_1=G, \qquad \theta_2 = K+\frac{4}{3}G. $$

The first M-eigenvalue $\theta_1$ is the shear modulus $G$ (or $\mu$), and %must be positive.
the second M-eigenvalue $\theta_2$ is the P-wave modulus $K+\frac{4}{3}G$. Let $\rho$ be the density, then $v_P :=\sqrt{\theta_2 /\rho} = \sqrt{(K+\frac{4}{3}G) /\rho}$ and $v_S := \sqrt{\theta_1 /\rho} = \sqrt{G /\rho}$ are the velocities of so-called longitudinal and shear elastic waves, and in the seismological literature, the corresponding plane waves are called P-waves and S-waves \cite[P.241]{Bower_book}.

Using the shear modulus $\mu$ and %the
Lam\'{e}'s first parameter $\lambda=K-\frac{2}{3} \mu $, the M-eigenvalues can be rewritten as
$$ \theta_1 = \mu, \qquad \theta_2 = \lambda + 2 \mu .$$
Although the shear modulus  $\mu$  must be positive,  Lam\'{e}'s first parameter $\lambda$ can be negative in principle; however, for most materials it is also positive.
Hence, in the isotropic case, strong ellipticity holds if and only if $\theta_2>0$, i.e., $\lambda > -2\mu$.

Note that for the isotropic case the elasticity tensor can be rewritten as
$C_{ijkl}=\lambda \delta_{ij} \delta_{kl}+\mu (\delta_{ik} \delta_{jl}+\delta_{il} \delta_{jk})$. Direct calculations yield that
$$\mathcal{C} xyxy = (\lambda + 2\mu) \langle x, y \rangle ^2 + \mu \left( \langle x, x \rangle \langle y, y \rangle - \langle x, y \rangle ^2 \right) .$$
From this formula we can also derive the strong ellipticity conditions: $\mu>0$ and $\lambda > -2\mu$.
Besides, we can derive that the Z-eigenvalue for this case is $\eta = K+\frac{4}{3}G= \lambda + 2 \mu$.
The positive definiteness of the strain energy yields $\mu >0$ and $\lambda + \frac{2}{3}\mu >0$ \cite[P.56]{Ting_Book96}. Obviously, the Z-eigenvalue gives a necessary condition, and the positive definiteness provides a sufficient condition for the strong ellipticity, respectively.

\section{Cubic system}

Cubic crystals, the simplest anisotropic case, can be described by three independent elasticity constants. In a properly chosen coordinate system, they can be put into the following Voigt matrix:
$$ ( c_{\alpha \beta } ) = {\begin{bmatrix}c_{11}&c_{12}&c_{12}&0&0&0\\c_{12}&c_{11}&c_{12}&0&0&0\\
c_{12}&c_{12}&c_{11}&0&0&0\\0&0&0&c_{66}&0&0\\0&0&0&0&c_{66}&0\\0&0&0&0&0&c_{66}\end{bmatrix}}.$$

Define $c_{11} := \beta_1$,  $c_{66} := \beta_2$ and $c_{12} := \beta_3$. The M-eigenvalue problem is given as follows.
\begin{eqnarray*}
x_1 ( \beta_1 y_1^2 + \beta_2 y_2^2 + \beta_2 y_3^2 ) + (\beta_3 + \beta_2) y_1 (x_2y_2 +  x_3y_3 ) &=& \theta x_1 , \\
x_2 ( \beta_2 y_1^2 + \beta_1 y_2^2 + \beta_2 y_3^2 ) + (\beta_3 + \beta_2) y_2 (x_1y_1 + x_3y_3 ) &=& \theta x_2 , \\
x_3 ( \beta_2 y_1^2 + \beta_2 y_2^2 + \beta_1 y_3^2 ) + (\beta_3 + \beta_2) y_3 (x_1y_1 + x_2y_2 ) &=& \theta x_3 , \\
y_1 ( \beta_1 x_1^2 + \beta_2 x_2^2 + \beta_2 x_3^2 ) + (\beta_3 + \beta_2) x_1 (x_2y_2 + x_3y_3 ) &=& \theta y_1 , \\
y_2 ( \beta_2 x_1^2 + \beta_1 x_2^2 + \beta_2 x_3^2 ) + (\beta_3 + \beta_2) x_2 (x_1y_1 + x_3y_3 ) &=& \theta y_2 , \\
y_3 ( \beta_2 x_1^2 + \beta_2 x_2^2 + \beta_1 x_3^2 ) +  (\beta_3 + \beta_2) x_3 (x_1y_1 + x_2y_2 ) &=& \theta y_3 .
\end{eqnarray*}

We can observe the symmetry in the polynomials above. The system of polynomial equations is unchanged if we perform the following operations: (1) Exchanging $x_i$ with $y_i$, denoted by $x_i \leftrightarrow y_i$ ($i=1,2,3$); (2) Exchanging $x_i$ with $x_j$ and meanwhile $y_i$ with $y_j$ ($x_i \leftrightarrow x_j$ and $y_i \leftrightarrow y_j$); (3) Cyclic permutations of the variables ($x_1 \rightarrow x_2 \rightarrow x_3 \rightarrow x_1$ and $y_1 \rightarrow y_2 \rightarrow y_3 \rightarrow y_1$).

Take into account the constraints $x_1^2+x_2^2+x_3^2 = 1$ and $y_1^2+y_2^2+y_3^2 = 1$. It is equivalent to
\begin{eqnarray*}
(\beta_1-\beta_2) x_1  y_1^2 + (\beta_3 + \beta_2) y_1  (x_2y_2 + x_3y_3) &=& (\theta - \beta_2) x_1 , \\
(\beta_1-\beta_2) x_2  y_2^2 + (\beta_3 + \beta_2) y_2  (x_1y_1 + x_3y_3) &=& (\theta - \beta_2) x_2 , \\
(\beta_1-\beta_2) x_3  y_3^2  + (\beta_3 + \beta_2) y_3 (x_1y_1 + x_2y_2 ) &=& (\theta - \beta_2) x_3 , \\
(\beta_1-\beta_2) y_1  x_1^2 +  (\beta_3 + \beta_2) x_1 (x_2y_2 + x_3y_3) &=& (\theta - \beta_2) y_1 , \\
(\beta_1-\beta_2) y_2  x_2^2 +  (\beta_3 + \beta_2) x_2 (x_1y_1 + x_3y_3) &=& (\theta - \beta_2) y_2 , \\
(\beta_1-\beta_2) y_3  x_3^2 +  (\beta_3 + \beta_2) x_3 (x_1y_1 + x_2y_2 ) &=& (\theta - \beta_2) y_3 .
\end{eqnarray*}
%That is,
%\begin{eqnarray*}
%y_1 [ (\beta_1-\beta_2) x_1  y_1 + (\beta_3 + \beta_2) (x_2y_2 + x_3y_3) ] &=& (\theta - \beta_2) x_1 , \\
%y_2 [ (\beta_1-\beta_2) x_2  y_2 + (\beta_3 + \beta_2) (x_1y_1 + x_3y_3) ] &=& (\theta - \beta_2) x_2 , \\
%y_3 [ (\beta_1-\beta_2) x_3  y_3 + (\beta_3 + \beta_2) (x_1y_1 + x_2y_2) ] &=& (\theta - \beta_2) x_3 ; \\
%x_1 [ (\beta_1-\beta_2) y_1  x_1 +  (\beta_3 + \beta_2) (x_2y_2 + x_3y_3) ] &=& (\theta - \beta_2) y_1 , \\
%x_2 [ (\beta_1-\beta_2) y_2  x_2 +  (\beta_3 + \beta_2) (x_1y_1 + x_3y_3) ] &=& (\theta - \beta_2) y_2 , \\
%x_3 [ (\beta_1-\beta_2) y_3  x_3 +  (\beta_3 + \beta_2) (x_1y_1 + x_2y_2) ]&=& (\theta - \beta_2) y_3 .
%\end{eqnarray*}
It is easy to verify that
$$ (\theta - \beta_2) x_i^2 = (\theta - \beta_2) y_i^2 = x_i y_i [ (\beta_1-\beta_2) x_i y_i +  (\beta_3 + \beta_2) (x_j y_j + x_k y_k) ],  $$
where $(i,j,k) = (1,2,3), (2,3,1)$, or $(3,1,2)$, and no summation convention used here for the subscripts $i,j$ and $k$. It is obvious that $\theta = \beta_2$, $x=(1,0,0)$ and $y=(0,1,0)$ are one admissible solution. What's more, we have $|x_i| = |y_i| ~ (i=1,2,3)$. In the following, we will consider different cases.

(I) Three components equal. That is, $x=y$. The system of polynomial equations reduce to
\begin{eqnarray*}
(\beta_1-\beta_2) x_1^3 + (\beta_3 + \beta_2) x_1 (x_2^2 + x_3^3)  &=& (\theta - \beta_2) x_1 , \\
(\beta_1-\beta_2) x_2^3 + (\beta_3 + \beta_2) x_2 (x_1^2 + x_3^2)  &=& (\theta - \beta_2) x_2 , \\
(\beta_1-\beta_2) x_3^3 + (\beta_3 + \beta_2) x_3 (x_1^2 + x_2^2)  &=& (\theta - \beta_2) x_3 .
\end{eqnarray*}
If $x_i~(i=1,2,3)$ are nonzeros, i.e., $x_1 x_2 x_3 \neq 0$, we have $\theta = \frac{\beta_1+4\beta_2+2\beta_3}{3}$. If only one component in $x$ is zero, then we obtain $\theta = \frac{\beta_1+2\beta_2 +\beta_3}{2}$. If two entries are zeros, then we get $\theta = \beta_1$.

(II) Two components equal. Considering the cyclic permutation symmetry of the polynomials, we set, for example, $x_1=y_1, x_2=y_2$ and $x_3=-y_3$. The system of polynomial equations become
\begin{eqnarray*}
x_1 [ (\beta_1-\beta_2) x_1^2 + (\beta_3 + \beta_2) (x_2^2 - x_3^3)  &=& (\theta - \beta_2) x_1 , \\
x_2 [ (\beta_1-\beta_2) x_2^2 + (\beta_3 + \beta_2) (x_1^2 - x_3^2)  &=& (\theta - \beta_2) x_2 , \\
(\beta_1 + \beta_3) x_3^2  &=& \theta + \beta_3 .
\end{eqnarray*}
Here we use $x_3 \neq 0$, otherwise it reduces the former case. If $x_1 x_2 \neq 0$, from the first two equations above, we have $ (\beta_1-\beta_2) (1-x_3^2) + (\beta_3 + \beta_2) (1 - 3 x_3^2)  = 2(\theta - \beta_2) $.
Combining with the third equation above, then we have
$$\theta = \frac{\beta_1^2 + 2\beta_1\beta_2 + \beta_1\beta_3  -2\beta_3^2}{3\beta_1+2\beta_2+5\beta_3},  $$
$$ x_3^2 = \frac{\beta_1+2\beta_2+3\beta_3}{3\beta_1+2\beta_2+5\beta_3}, \quad x_2^2 = x_1^2 =  \frac{\beta_1+\beta_3}{3\beta_1+2\beta_2+5\beta_3}.$$
If $x_2 = 0, x_1 \neq 0$, we have $\theta = \beta_1$; If $x_1 = 0, x_2 \neq 0$, we have $\theta =  (\beta_1-\beta_3)/2 $.

(III) One component equals. For example, $x_1=y_1, x_2=-y_2$ and $x_3=-y_3$. The system of polynomial equations read
\begin{eqnarray*}
(\beta_1 + \beta_3) x_1^3  &=& (\theta + \beta_3) x_1, \\
x_2 [ (\beta_1-\beta_2) x_2^2 - (\beta_3 + \beta_2) (x_1^2 - x_3^2)  &=& (\theta - \beta_2) x_2 ,\\
x_3 [ (\beta_1-\beta_2) x_3^2 - (\beta_3 + \beta_2) (x_1^2 - x_2^3)  &=& (\theta - \beta_2) x_3 .
\end{eqnarray*}
For $x_1 \neq 0$, we find that $\theta = \frac{\beta_1^2 + 2\beta_1\beta_2 + \beta_1\beta_3  -2\beta_3^2}{3\beta_1+2\beta_2+5\beta_3}$ when $x_2 x_3 \neq 0$, and $\theta = \frac{\beta_1-\beta_3}{2}$ when $x_2 x_3 = 0$. For $x_1 = 0$, we find that $\theta = \frac{\beta_1+2\beta_2 +\beta_3}{2}$ when $x_2 x_3 \neq 0$, and $\theta = \beta_1$ when $x_2 x_3 = 0$.
This case does not yield any new M-eigenvalues.
By the way, the case where $x=-y$ can be reduced to the first case, and does not give new M-eigenvalues either.

In all, we can find the following M-eigenvalues:
%\footnote{For the special case where $x=y$, we have $\theta=(\beta_1+2\beta_3+4\beta_2)/3$.
%\\ The last two eigenvalues will give us new information? To be checked later... }
\begin{eqnarray*}
\theta_1 &=& \beta_1, \\
\theta_2 &=& \beta_2, \\
\theta_3 &=& \frac{\beta_1-\beta_3}{2}, \\
\theta_4 &=& \frac{\beta_1+2\beta_2 +\beta_3}{2} , \\
\theta_5 &=& \frac{\beta_1+4\beta_2+2\beta_3}{3}, \\
\theta_6 &=& \frac{\beta_1^2 + 2\beta_1\beta_2 + \beta_1\beta_3  -2\beta_3^2}{3\beta_1+2\beta_2+5\beta_3}.
\end{eqnarray*}

For the cubic system the Z-eigenvalue problem \eqref{eqn:ZeigDef} reduces to
$$(\beta_1-2\beta_2-\beta_3) x_j^3 = (\eta -2\beta_2-\beta_3) x_j \quad (j=1,2,3) . $$
Hence, we have the Z-eigenvalues: $\eta_1 = \beta_1$, $\eta_2 = \frac{\beta_1+\beta_3}{2} +\beta_2$, and $\eta_3 = \frac{\beta_1+2\beta_3+4\beta_2}{3}$.

Let $\nu$ be the Poisson's ratio, and $E$ the Young's modulus.  Then $\beta_1=c_{11} = E(1-\nu)/(1-\nu-2\nu^2)$, $\beta_2 = c_{66} = \mu %= \frac{E}{2(1+\nu)}
$, and $\beta_3 = c_{12} = E \nu/(1-\nu-2\nu^2)$.   The ratio
$$ A:= \frac{\theta_2}{\theta_3} = \frac{2 \beta_2}{\beta_1 - \beta_3} = \frac{2 c_{66}}{c_{11}-c_{12}} = \frac{2 \mu (1+\nu)}{E} $$
provides a convenient measure of anisotropy \cite[p.87]{Bower_book}, often referred to as the anisotropy factor. For $A = 1$, the material is isotropic.
The M-eigenvalues $\theta_1$, $\theta_2$ are obviously the elasticity constants, and the M-eigenvalues $\theta_3$, $\theta_4$, and $\theta_5$ are related to distinct wave types \cite[p.269]{Haussuhl_book}.

For the strong ellipticity of cubic material, the necessary %and sufficient
conditions are that $\theta_i > 0 ~(i=1,\cdots,5)$, i.e.,
\begin{equation}\label{eqn:Cubic1}
c_{11} > 0, \quad c_{66} > 0, \quad c_{11} > c_{12},
\end{equation}
\begin{equation}\label{eqn:Cubic2}
c_{11}+2c_{66}+c_{12} > 0, \quad c_{11} + 4 c_{66} + 2 c_{12} > 0.
\end{equation}

Suppose that the absolute value $\beta_3$ is small, then $\theta_6 \approx \frac{2\beta_1 \theta_4}{3\beta_1 + 2\beta_2} >0$. Note that the numerator of $\theta_6$ is $2\beta_1\beta_2 + (\beta_1-\beta_3)(\beta_1+2\beta_3)$ and the denominator of $\theta_6$ is $\frac{1}{2}\beta_1 + 2\beta_2 + \frac{5}{2}(\beta_1+2\beta_3)$. Thus, $\beta_1+2\beta_3>0$ is sufficient to guarantee $\theta_6>0$. Using the elasticity constants, this condition is reexpressed as
\begin{equation}\label{eqn:Cubic3}
c_{11} + 2 c_{12} > 0.
\end{equation}
The quantity $c_{11} + 2 c_{12}$ is the bulk modulus for three-dimensional deformations, and the condition above is necessary for the positive definiteness of elasticity tensor \cite[P.51]{PaoloVannucci_book18} \cite[P.59]{Ting_Book96}.
We note that the conditions \eqref{eqn:Cubic1}, \eqref{eqn:Cubic2} and \eqref{eqn:Cubic3} are sufficient for the strong ellipticity. The cubic materials, such as Au, Cu, NaCl, MgAl$_2$O$_4$, FeO, MgO, FeS$_2$, Fe$_3$O$_4$, etc., satisfy these conditions \cite[P.46]{Ahrens_Book1995}.

{\bf Remarks.}
Using the positiveness of the M-eigenvalues $\theta_3$ and $\theta_4$, we have $4\theta_3\theta_4 =(\beta_1-\beta_3)(\beta_1+\beta_3+2\beta_2)>0$, that is, $(\beta_3+\beta_2)^2 < (\beta_1+\beta_2)^2$.
Hence from the positiveness of first four M-eigenvalues, we have the necessary conditions
$$\beta_1>0, \quad \beta_2>0, \quad |\beta_3 + \beta_2| < \beta_1 + \beta_2. $$
These are equivalent to the already known results, the strong ellipticity conditions for cubic system \cite{SChirita_JE07}: 
$$c_{11}>0, \quad c_{66}>0, \quad |c_{12} + c_{66}| < c_{11} + c_{66}. $$

\section{Polar anisotropy}

The polar anisotropy is the best known anisotropic case, also called the transverse isotropy.
Many materials belong to this class, for instance, wood, fiber reinforced composites, laminated steel and so on.
The corresponding Voigt matrix has the following form:
$$ ( c_{\alpha \beta } ) = {\begin{bmatrix}c_{11} &c_{11}-2c_{66} & c_{13}&0&0&0 \\c_{11}-2c_{66}&c_{11}&c_{13}&0&0&0 \\ c_{13}&c_{13}&c_{33}&0&0&0 \\0&0&0&c_{44}&0&0\\0&0&0&0&c_{44}&0\\0&0&0&0&0&c_{66}\end{bmatrix}}.$$
The elasticity tensor is described by five transversely isotropic parameters:
$c_{11}=c_{22} = C_{1111} =C_{2222} := \alpha_1$, $c_{33}=C_{3333} := \alpha_2$, $c_{44}=c_{55}=C_{2323}=C_{3131} :=\alpha_3$, $c_{66}=C_{1212} :=\alpha_4$, $c_{23}=c_{13}=C_{2233}=C_{1133} := \alpha_5$, $c_{12}=C_{1122}=\alpha_1-2\alpha_4$.

The corresponding M-eigenvalue problem reads as follows.
\begin{eqnarray*}
x_1 ( \alpha_1 y_1^2 + \alpha_4 y_2^2 + \alpha_3 y_3^2 ) + y_1 [ (\alpha_1 -\alpha_4)x_2y_2 + (\alpha_3+\alpha_5)x_3y_3 ] &=& \theta x_1 , \\
x_2 ( \alpha_1 y_2^2 + \alpha_4 y_1^2 + \alpha_3 y_3^2 ) + y_2 [ (\alpha_1 -\alpha_4)x_1y_1 + (\alpha_3+\alpha_5)x_3y_3 ] &=& \theta x_2 , \\
x_3 ( \alpha_3 y_1^2 + \alpha_3 y_2^2 + \alpha_2 y_3^2 ) + (\alpha_3+\alpha_5) y_3 (x_1y_1 + x_2y_2 ) &=& \theta x_3 , \\
y_1 ( \alpha_1 x_1^2 + \alpha_4 x_2^2 + \alpha_3 x_3^2 ) + x_1 [ (\alpha_1 -\alpha_4)x_2y_2 + (\alpha_3+\alpha_5)x_3y_3 ] &=& \theta y_1 , \\
y_2 (  \alpha_1 x_2^2 + \alpha_4 x_1^2 + \alpha_3 x_3^2 ) + x_2 [ (\alpha_1 -\alpha_4)x_1y_1 + (\alpha_3+\alpha_5)x_3y_3 ] &=& \theta y_2 , \\
y_3 ( \alpha_3 x_1^2 + \alpha_3 x_2^2 + \alpha_2 x_3^2 ) + (\alpha_3+\alpha_5) x_3 (x_1y_1 + x_2y_2 ) &=& \theta y_3 .
\end{eqnarray*}

Take into account the constraints $x_1^2+x_2^2+x_3^2 = 1$ and $y_1^2+y_2^2+y_3^2 = 1$. It is equivalent to
\begin{eqnarray*}
x_1 [ (\alpha_1-\alpha_3) y_1^2 + (\alpha_4-\alpha_3) y_2^2 ] + y_1 [ (\alpha_1 -\alpha_4)x_2y_2 + (\alpha_3+\alpha_5)x_3y_3 ] &=& (\theta - \alpha_3) x_1 , \\
x_2 [ (\alpha_1-\alpha_3) y_2^2 + (\alpha_4-\alpha_3) y_1^2 ] + y_2 [ (\alpha_1 -\alpha_4)x_1y_1 + (\alpha_3+\alpha_5)x_3y_3 ] &=& (\theta - \alpha_3) x_2 , \\
(\alpha_2-\alpha_3) x_3  y_3^2  + (\alpha_3+\alpha_5) y_3 (x_1y_1 + x_2y_2 ) &=& (\theta - \alpha_3) x_3 , \\
y_1 [ (\alpha_1-\alpha_3) x_1^2 + (\alpha_4-\alpha_3) x_2^2 ]  + x_1 [ (\alpha_1 -\alpha_4)x_2y_2 + (\alpha_3+\alpha_5)x_3y_3 ] &=& (\theta - \alpha_3) y_1 , \\
y_2 [ (\alpha_1-\alpha_3) x_2^2 + (\alpha_4-\alpha_3) x_1^2 ] + x_2 [ (\alpha_1 -\alpha_4)x_1y_1 + (\alpha_3+\alpha_5)x_3y_3 ] &=& (\theta - \alpha_3) y_2 , \\
(\alpha_2-\alpha_3) y_3  x_3^2   + (\alpha_3+\alpha_5) x_3 (x_1y_1 + x_2y_2 ) &=& (\theta - \alpha_3) y_3 .
\end{eqnarray*}
We can verify that
  $$(\theta - \alpha_3) x_3^2 = (\theta - \alpha_3) y_3^2 = (\alpha_2-\alpha_3) x_3^2 y_3^2   + (\alpha_3+\alpha_5) x_3 y_3 (x_1y_1 + x_2y_2 ). $$
Hence, $\theta = \alpha_3$, or $ |x_3| = |y_3|$.  For the case where $x_3 = y_3 =0$, we obtain $\theta = \alpha_1$ and $\alpha_4$. Then in the following we only consider $ |x_3| = |y_3| \neq 0$.

(I) The case where $x_3 = y_3$. We need to solve 
\begin{eqnarray*}
x_1 [ (\alpha_1-\alpha_3) y_1^2 + (\alpha_4-\alpha_3) y_2^2 ] + y_1 [ (\alpha_1 -\alpha_4)x_2y_2 + (\alpha_3+\alpha_5)x_3^2 ] &=& (\theta - \alpha_3) x_1 , \\
x_2 [ (\alpha_1-\alpha_3) y_2^2 + (\alpha_4-\alpha_3) y_1^2 ] + y_2 [ (\alpha_1 -\alpha_4)x_1y_1 + (\alpha_3+\alpha_5)x_3^2 ] &=& (\theta - \alpha_3) x_2 , \\
y_1 [ (\alpha_1-\alpha_3) x_1^2 + (\alpha_4-\alpha_3) x_2^2 ]  + x_1 [ (\alpha_1 -\alpha_4)x_2y_2 + (\alpha_3+\alpha_5)x_3^2 ] &=& (\theta - \alpha_3) y_1 , \\
y_2 [ (\alpha_1-\alpha_3) x_2^2 + (\alpha_4-\alpha_3) x_1^2 ] + x_2 [ (\alpha_1 -\alpha_4)x_1y_1 + (\alpha_3+\alpha_5)x_3^2 ] &=& (\theta - \alpha_3) y_2 , \\
(\alpha_2-\alpha_3) x_3^3  + (\alpha_3+\alpha_5) x_3 (x_1y_1 + x_2y_2 ) &=& (\theta - \alpha_3) x_3 .
\end{eqnarray*}
We observe that the system of polynomial equations is symmetric when $x_1$ exchanges with $x_2$ and $y_1$ exchanges with $y_2$ ($x_1 \leftrightarrow x_2$ and $y_1 \leftrightarrow y_2$).
We therefore set $x_1 = x_2$ and $y_1 = y_2$, and the unknowns are $x_1, y_1, x_3$ and $\theta$. Then we have
\begin{eqnarray*}
(\alpha_1-2\alpha_3+\alpha_4) x_1 y_1^2 + y_1 [ (\alpha_1 -\alpha_4)x_1y_1 + (\alpha_3+\alpha_5)x_3^2 ] &=& (\theta - \alpha_3) x_1 , \\
(\alpha_1-2\alpha_3+\alpha_4) y_1 x_1^2 + x_1 [ (\alpha_1 -\alpha_4)x_1y_1 + (\alpha_3+\alpha_5)x_3^2 ] &=& (\theta - \alpha_3) y_1 , \\
(\alpha_2-\alpha_3) x_3^3  + 2(\alpha_3+\alpha_5) x_3 x_1y_1  &=& (\theta - \alpha_3) x_3 .
\end{eqnarray*}
We have the M-eigenvalues: $\theta = \alpha_1$, $\theta = \alpha_2$, $\theta = \frac{\alpha_1\alpha_2-\alpha_5^2}{\alpha_1+\alpha_2+2\alpha_5} $, $\theta = \frac{\alpha_1\alpha_2 - (2\alpha_3+\alpha_5)^2}{\alpha_1+\alpha_2 - 2(2\alpha_3+\alpha_5)}$.

Using the symmetry denoted by $x_1 \leftrightarrow y_1$ and $x_2 \leftrightarrow y_2$, we have the system of polynomial equations with unknowns $x_1, x_2, x_3$ and $\theta$.
\begin{eqnarray*}
x_1 [ (\alpha_1-\alpha_3) x_1^2 + (\alpha_4-\alpha_3) x_2^2 ] + x_1 [ (\alpha_1 -\alpha_4)x_2^2 + (\alpha_3+\alpha_5)x_3^2 ] &=& (\theta - \alpha_3) x_1 , \\
x_2 [ (\alpha_1-\alpha_3) x_2^2 + (\alpha_4-\alpha_3) x_1^2 ] + x_2 [ (\alpha_1 -\alpha_4)x_1^2 + (\alpha_3+\alpha_5)x_3^2 ] &=& (\theta - \alpha_3) x_2 , \\
(\alpha_2-\alpha_3) x_3^3  + (\alpha_3+\alpha_5) x_3 (x_1^2 + x_2^2 ) &=& (\theta - \alpha_3) x_3 .
\end{eqnarray*}
We have the M-eigenvalues:
$\theta = \alpha_1$, $\theta = \alpha_2$, $\theta = \frac{\alpha_1\alpha_2 - (2\alpha_3+\alpha_5)^2}{\alpha_1+\alpha_2 - 2(2\alpha_3+\alpha_5)}$.

\footnote{
The system of polynomial equations is symmetric about $x_1$ and $\pm y_1$. In fact, the conditions that $x_1 = x_2$ and $y_1 = y_2$ also imply that $|x_1| = |y_1|$. Then the unknowns are $x_1, x_3$ and $\theta$, and the system of polynomial equations read
\begin{eqnarray*}
2(\alpha_1-\alpha_3)  x_1^3 + (\alpha_3+\alpha_5)x_1x_3^2 &=& (\theta - \alpha_3) x_1 , \\
(\alpha_2-2\alpha_3-\alpha_5) x_3^3   &=& (\theta - 2\alpha_3 - \alpha_5) x_3 ,
\end{eqnarray*}
together with the constraint $2 x_1^2 + x_3^2 = 1$. When $x_1 x_3 \neq 0$, we have
$$\theta = \frac{\alpha_1\alpha_2 - (2\alpha_3+\alpha_5)^2}{\alpha_1+\alpha_2 - 2(2\alpha_3+\alpha_5)}.$$
%$$ \theta = \frac{\alpha_3(\alpha_1-2\alpha_3-\alpha_5) + \alpha_1(\alpha_2-2\alpha_3-\alpha_5)}{\alpha_1+\alpha_2 - 2(2\alpha_3+\alpha_5)}. $$
%
When $x_1=0$, that is $x_3=1$, then we have $\theta = \alpha_2$. In addition, we have $\theta=\alpha_1$.
}

(II) The case where $x_3 = - y_3$. We need to solve
\begin{eqnarray*}
x_1 [ (\alpha_1-\alpha_3) y_1^2 + (\alpha_4-\alpha_3) y_2^2 ] + y_1 [ (\alpha_1 -\alpha_4)x_2y_2 - (\alpha_3+\alpha_5)x_3^2 ] &=& (\theta - \alpha_3) x_1 , \\
x_2 [ (\alpha_1-\alpha_3) y_2^2 + (\alpha_4-\alpha_3) y_1^2 ] + y_2 [ (\alpha_1 -\alpha_4)x_1y_1 - (\alpha_3+\alpha_5)x_3^2 ] &=& (\theta - \alpha_3) x_2 , \\
y_1 [ (\alpha_1-\alpha_3) x_1^2 + (\alpha_4-\alpha_3) x_2^2 ]  + x_1 [ (\alpha_1 -\alpha_4)x_2y_2 - (\alpha_3+\alpha_5)x_3^2 ] &=& (\theta - \alpha_3) y_1 , \\
y_2 [ (\alpha_1-\alpha_3) x_2^2 + (\alpha_4-\alpha_3) x_1^2 ] + x_2 [ (\alpha_1 -\alpha_4)x_1y_1 - (\alpha_3+\alpha_5)x_3^2 ] &=& (\theta - \alpha_3) y_2 , \\
(\alpha_2-\alpha_3) x_3^3  - (\alpha_3+\alpha_5) x_3 (x_1y_1 + x_2y_2 ) &=& (\theta - \alpha_3) x_3 .
\end{eqnarray*}

Using the symmetry denoted by $x_1 \leftrightarrow x_2$ and $y_1 \leftrightarrow y_2$, we have
\begin{eqnarray*}
(\alpha_1-2\alpha_3+\alpha_4) x_1 y_1^2  + y_1 [ (\alpha_1 -\alpha_4)x_1y_1 - (\alpha_3+\alpha_5)x_3^2 ] &=& (\theta - \alpha_3) x_1 , \\
(\alpha_1-2\alpha_3+\alpha_4)y_1 x_1^2   + x_1 [ (\alpha_1 -\alpha_4)x_1y_1 - (\alpha_3+\alpha_5)x_3^2 ] &=& (\theta - \alpha_3) y_1 , \\
(\alpha_2-\alpha_3) x_3^3  - 2(\alpha_3+\alpha_5) x_3 x_1y_1  &=& (\theta - \alpha_3) x_3 .
\end{eqnarray*}
We have the M-eigenvalues:
$\theta=\alpha_1$, $\theta=\alpha_2$, $\theta = \frac{\alpha_1\alpha_2-\alpha_5^2}{\alpha_1+\alpha_2+2\alpha_5} $, $\theta = \frac{\alpha_1\alpha_2 - (2\alpha_3+\alpha_5)^2}{\alpha_1+\alpha_2 - 2(2\alpha_3+\alpha_5)}$.

Using the symmetry denoted by $x_1 \leftrightarrow y_1$ and $x_2 \leftrightarrow y_2$, we have
\begin{eqnarray*}
x_1 [ (\alpha_1-\alpha_3) x_1^2 + (\alpha_4-\alpha_3) x_2^2 ] + x_1 [ (\alpha_1 -\alpha_4)x_2^2 - (\alpha_3+\alpha_5)x_3^2 ] &=& (\theta - \alpha_3) x_1 , \\
x_2 [ (\alpha_1-\alpha_3) x_2^2 + (\alpha_4-\alpha_3) x_1^2 ] + x_2 [ (\alpha_1 -\alpha_4)x_1^2 - (\alpha_3+\alpha_5)x_3^2 ] &=& (\theta - \alpha_3) x_2 , \\
(\alpha_2-\alpha_3) x_3^3  - (\alpha_3+\alpha_5) x_3 (x_1^2 + x_2^2 ) &=& (\theta - \alpha_3) x_3 .
\end{eqnarray*}
We have the M-eigenvalues:
$\theta=\alpha_1$, $\theta=\alpha_2$, $\theta = \frac{\alpha_1\alpha_2-\alpha_5^2}{\alpha_1+\alpha_2+2\alpha_5} $.

\footnote{
Observing that the system of polynomial equations is symmetric when $x_1$ exchanges with $x_2$ and $y_1$ exchanges with $y_2$ ($x_1 \leftrightarrow x_2$ and $y_1 \leftrightarrow y_2$),
and noticing that
$x_1 = x_2$ and $y_1 = y_2$ also imply that $x_1=y_1$, we have the system about three unknowns  $x_1, x_3$ and $\theta$.
\begin{eqnarray*}
2(\alpha_1-\alpha_3)  x_1^3 - (\alpha_3+\alpha_5)x_1x_3^2 &=& (\theta - \alpha_3) x_1 , \\
(\alpha_2+\alpha_5) x_3^3   &=& (\theta + \alpha_5) x_3.
\end{eqnarray*}
Therefore, we have the M-eigenvalues: $\theta = \alpha_2$ and $\theta= \frac{\alpha_1\alpha_2-\alpha_5^2}{\alpha_1+\alpha_2+2\alpha_5} $. In addition, we have $\theta=\alpha_1$.
}

%------------------------------------------------------

%-*** symmetric polynomials ***-
Define the following notations:
$$\Lambda= \mathrm{diag}(1,0,-1), \quad D=\mathrm{diag}(0,1,-1), \quad P = \begin{pmatrix} 0 & 1 & 0 \\ 0 & 0 & 1 \\ 1 & 0 & 0 \end{pmatrix}; $$
\begin{eqnarray*}
\zeta= \sqrt{ \frac{2 \alpha_3 + \alpha_5 - \alpha_1}{4 \alpha_3 + 2 \alpha_5 - \alpha_1 - \alpha_2} }, && \phi =  \sqrt{ \frac{\alpha_1 + \alpha_5}{ \alpha_1 + \alpha_2 + 2 \alpha_5} } , \\
\eta = \sqrt{ \frac{2 \alpha_3 + \alpha_5 - \alpha_2}{4 \alpha_3 + 2 \alpha_5 - \alpha_1 - \alpha_2} }, && \psi =  \sqrt{ \frac{\alpha_2 + \alpha_5}{ \alpha_1 + \alpha_2 + 2 \alpha_5} } .
%
%\theta_*=\frac{\alpha_1\alpha_2 - (2\alpha_3+\alpha_5)^2}{\alpha_1+\alpha_2 - 2(2\alpha_3+\alpha_5)} , && \tau_* = \frac{\alpha_1\alpha_2-\alpha_5^2}{\alpha_1+\alpha_2+2\alpha_5}.
\end{eqnarray*}
We list 62 M-eigenvalue pairs and show them in Table \ref{tab:Meig4PolarAniso}.
We find the following six M-eigenvalues $\theta_i$:
%$$\alpha_1, \quad \alpha_2, \quad \alpha_3, \quad \alpha_4, \quad \theta_*, \quad \tau_* .$$
\begin{eqnarray*}
%\theta_1 &=& \alpha_1, \\
%\theta_2 &=& \alpha_2, \\
%\theta_3 &=& \alpha_3, \\
%\theta_4 &=& \alpha_4, \\
\theta_i &=& \alpha_i, \quad (i = 1,2,3,4) \\
\theta_5 &=& \frac{\alpha_1\alpha_2-\alpha_5^2}{\alpha_1+\alpha_2+2\alpha_5} , \\
\theta_6 &=& \frac{\alpha_1\alpha_2 - (2\alpha_3+\alpha_5)^2}{\alpha_1+\alpha_2 - 2(2\alpha_3+\alpha_5)}.
\end{eqnarray*}

\newsavebox{\tableboxMeig}
\begin{lrbox}{\tableboxMeig}

\begin{tabular}{cccccccccccccc}\hline
$x$ &
$\begin{pmatrix}\pm 1\\0\\0\end{pmatrix}$ & $\begin{pmatrix}0\\\pm 1\\0\end{pmatrix}$ & $\begin{pmatrix}0\\0\\\pm 1\end{pmatrix}$ & $\begin{pmatrix}0\\0\\\pm 1\end{pmatrix}$   & $\begin{pmatrix}\pm 1\\0\\0\end{pmatrix}$ &
$\begin{pmatrix}0\\ 0\\\pm 1\end{pmatrix}$  & $\begin{pmatrix}0\\ \pm 1\\0\end{pmatrix}$ &
$\begin{pmatrix}0\\ \pm 1\\0\end{pmatrix}$  & $\begin{pmatrix}\pm 1\\ 0\\0\end{pmatrix}$ &
$\begin{pmatrix}\pm \psi\\ 0 \\ \pm \phi \end{pmatrix}$ &
$\begin{pmatrix}0\\\pm \psi \\ \pm \phi \end{pmatrix}$  &
$\begin{pmatrix}0\\\pm \eta \\ \pm \zeta \end{pmatrix}$ &
$\begin{pmatrix}\pm \eta \\ 0 \\ \pm \zeta \end{pmatrix}$
\\
$y$ & $x$ & $x$ & $x$ &
$\pm P^2 x$   & $\pm P x$ & $\pm P x$   & $\pm P^2 x$ & $\pm P x$   & $\pm P^2 x$ & $ \pm \Lambda x$  & $ \pm Dx$ & $ \pm x $  & $ \pm x $2
\\ \hline
$\theta$ & $\alpha_1$ & $\alpha_1$ & $\alpha_2$ & $\alpha_3$ & $\alpha_3$ & $\alpha_3$ & $\alpha_3$ & $\alpha_4$ & $\alpha_4$ & $\theta_5$  & $\theta_5$ & $\theta_6$ & $\theta_6$    \\ \hline
\end{tabular}

\end{lrbox}

\begin{table}
 \centering \scalebox{0.6}{\usebox{\tableboxMeig}}
 \caption{\label{tab:Meig4PolarAniso} M-eigenvalues and eigenvectors of polar anisotropic case.  }
\end{table}

The corresponding Z-eigenvalue problem reduces to
\begin{eqnarray*}
(-\alpha_1+2\alpha_3+ \alpha_5) x_1 x_3^2  &=& x_1 (\eta -\alpha_1)  , \\
(-\alpha_1+2\alpha_3+ \alpha_5) x_2 x_3^2  &=& x_2 (\eta -\alpha_1)  , \\
(\alpha_2 -2\alpha_3 - \alpha_5)  x_3^3  &=& x_3 (\eta -2\alpha_3 - \alpha_5)   .
\end{eqnarray*}
We can deduce the Z-eigenvalues: $\eta_1 = \alpha_1$, $\eta_2 = \alpha_2$ and $\eta_3 = \theta_6$.

Note that  $\alpha_k (k=1,\cdots,4)$ are the elasticity constants defined in terms of generalized Young's moduli, shear moduli and Poisson's ratios \cite[p.84]{Bower_book}.
That is, $\theta_1=\alpha_1=E_p(1-\nu_{pt}\nu_{tp}) \gamma$, $\theta_2=\alpha_2 = E_t(1-\nu_p^2)   \gamma$, $\theta_3=\alpha_3 = \mu_t$, $\theta_4=\alpha_4=\mu_p$ with $\gamma= 1/ (1-\nu_p^2 -2\nu_{tp}\nu_{tp} -2 \nu_p \nu_{pt} \nu_{tp}) $, where $E_p = E_1 = E_2$, $E_t=E_3$, $\nu_p = \nu_{12}=\nu_{21}$, $\nu_{tp}= \nu_{31} =\nu_{32}$, $\nu_{pt}=\nu_{13}=\nu_{23}$, $\mu_p=E_p/2(1+\nu_p)$, and
$E_k$ is the Young's modulus along axis $k$,  $\nu_{ij}$ is the Poisson's ratio  corresponding to a contraction in direction $j$ while an extension %is applied
in direction $i$. In addition, the elasticity constant $\alpha_5 = E_p (\nu_{tp}+\nu_p \nu_{tp}) \gamma$.

%-----$\theta_5=?,\theta_6=?$-------------------

Since the strong ellipticity is equivalent to the positiveness of the M-eigenvalues \cite{QDH09}, We now in the position to investigate the positiveness of M-eigenvalues. We need the basic fact: If $a>0, b>0$ then $\frac{ab-c^2}{a+b-2c}>0$ is equivalent to the conditions that $|c| < \sqrt{ab}$ or $ c > \frac{a+b}{2}$.
For $\theta_5 >0$, we have $| \alpha_5 | < \sqrt{\alpha_1 \alpha_2}$ or $\alpha_5 < - \frac{\alpha_1 + \alpha_2}{2}$.
For $\theta_6 >0$, we have $|2\alpha_3 + \alpha_5 | < \sqrt{\alpha_1 \alpha_2}$ or $2\alpha_3 + \alpha_5 > \frac{\alpha_1 + \alpha_2}{2}$.
All the Young's and shear moduli are strictly positive quantities, and usually Young's moduli are large while the shear moduli are small and Poisson's ratios are even smaller \cite{CowinYangMehrabadi_JE99,ItskovAksel_AM02,Lempriere_AIAA68} \cite[p.73, Table 3.1]{Bower_book} \cite[p.75]{Bower_book} \cite[p.285]{Haussuhl_book}.
Due to the practical considerations, we reject the conditions that $\alpha_5 < - \frac{\alpha_1 + \alpha_2}{2}$, % and $2\alpha_3 + \alpha_5 > \frac{\alpha_1 + \alpha_2}{2}$,
and provide the following %sufficient
conditions for the strong ellipticity.
$$ \alpha_i>0 ~(i=1,2,3,4), \quad | \alpha_5 | < \sqrt{\alpha_1 \alpha_2},$$
$$ |2\alpha_3 + \alpha_5 | < \sqrt{\alpha_1 \alpha_2}
\quad \text{or} \quad 2\alpha_3 + \alpha_5 > \frac{\alpha_1 + \alpha_2}{2}  .$$
Notice that $\alpha_1=c_{11}$, $\alpha_2 = c_{33}$, $\alpha_3=c_{55}$, $\alpha_5 = c_{13}$ and $\alpha_1-2\alpha_4 = c_{12}$. Therefore, we have the following %sufficient
conditions for strong ellipticity.
\begin{equation}\label{eqn:polar1}
 c_{11} > 0, \quad c_{33} > 0, \quad c_{55} > 0, \quad c_{11} > c_{12},
\end{equation}
\begin{equation}\label{eqn:polar2}
| c_{13}|   < \sqrt{c_{11} c_{33} }.
\end{equation}
\begin{equation}\label{eqn:polar3}
|2 c_{55} + c_{13}|  < \sqrt{c_{11} c_{33}} \quad \text{or} \quad 2 c_{55} + c_{13} > \frac{c_{11} + c_{33}}{2} .
\end{equation}
%
%Define $\bar{\alpha} = \alpha_3 + \alpha_5$. We say more about the condition that $\theta_5>0$:
%\begin{enumerate}[(i)]
%\item If $ \bar{\alpha}+\alpha_3 > \frac{\alpha_1+\alpha_2}{2}$, then obviously $\theta_5 >0$.
%
%\item If $ \bar{\alpha}+\alpha_3 < \frac{\alpha_1+\alpha_2}{2}$, then $\theta_5 >0$ is equivalent to  $(\bar{\alpha}+\alpha_3)^2 < \alpha_1 \alpha_2$. For this case, we obtain that
%$$ |\bar{\alpha}| - \alpha_3 \leq \left| ~ |\bar{\alpha}| - \alpha_3 ~\right| \leq |\bar{\alpha} + \alpha_3| < \sqrt{\alpha_1 \alpha_2}. $$
%
%\end{enumerate}
%
For example, the elastic moduli of  %a transversely isotropic monocrystalline zinc
Al$_2$O$_3$, TiB$_2$ and Zinc
satisfy the second inequality in \eqref{eqn:polar3} \cite[P.22]{DingChenZhang_Book2006} \cite{LubardaChen_JMMS08}; and the elastic constants of
GaS, GaSe, MoS$_2$, NbSe$_2$, SiC, Ti and ZnO satisfy the first condition in in \eqref{eqn:polar3} \cite[P.22]{DingChenZhang_Book2006}.

{\bf Remarks.}
Since $ |\alpha_3+\alpha_5| - \alpha_3 \leq \left| ~ |\alpha_3+\alpha_5| - \alpha_3 ~\right| \leq | \alpha_3+\alpha_5 + \alpha_3| < \sqrt{\alpha_1 \alpha_2}$. We therefore have $|\alpha_5 + \alpha_3| < \alpha_3 + \sqrt{\alpha_1 \alpha_2}$, i.e.,
\begin{equation}\label{eqn:Chirita}
|c_{13}+c_{55}|<c_{55}+\sqrt{c_{11}c_{33}}.
\end{equation}
Note that the inequalities \eqref{eqn:polar1} and \eqref{eqn:Chirita} already appeared in \cite{SChirita_JE07} as the strong ellipticity conditions of polar anisotropic case.

\section{Tetragonal}

With respect to a proper orthonormal basis, the tensor of tetragonal symmetry has seven nonzero elasticities. But by means of a coordinate transformation, it can be reduced to an elasticity tensor of tetragonal symmetry with six independent components.
The corresponding Voigt matrix has the following form:
$$ ( c_{\alpha \beta } ) = {\begin{bmatrix}c_{11} &c_{12} & c_{13}&0&0&0 \\c_{12}&c_{11}&c_{13}&0&0&0 \\ c_{13}&c_{13}&c_{33}&0&0&0 \\0&0&0&c_{44}&0&0\\0&0&0&0&c_{44}&0\\0&0&0&0&0&c_{66}\end{bmatrix}}.$$
The elasticity tensor is described by six parameters:
$c_{11}=c_{22}= C_{1111} = C_{2222}:=\alpha_1$,  $c_{33}=C_{3333} := \alpha_2$, $c_{44}=c_{55}=C_{2323}=C_{3232} =C_{1313}=C_{3131}:=\alpha_3$, $c_{66}=C_{1212}=C_{2121}:=\alpha_4$, $c_{13}=c_{23}=C_{2233} =C_{1133} := \alpha_5$, $c_{12}=C_{1122}:=\alpha_6$.

The corresponding M-eigenvalue problem reads as follows.
\begin{eqnarray*}
x_1 ( \alpha_1 y_1^2 + \alpha_4 y_2^2 + \alpha_3 y_3^2 ) + y_1 [ (\alpha_6 +\alpha_4)x_2y_2 + (\alpha_3+\alpha_5)x_3y_3 ] &=& \theta x_1 , \\
x_2 ( \alpha_1 y_2^2 + \alpha_4 y_1^2 + \alpha_3 y_3^2 ) + y_2 [ (\alpha_6 +\alpha_4)x_1y_1 + (\alpha_3+\alpha_5)x_3y_3 ] &=& \theta x_2 , \\
x_3 ( \alpha_3 y_1^2 + \alpha_3 y_2^2 + \alpha_2 y_3^2 ) + (\alpha_3+\alpha_5) y_3 (x_1y_1 + x_2y_2 ) &=& \theta x_3 , \\
y_1 ( \alpha_1 x_1^2 + \alpha_4 x_2^2 + \alpha_3 x_3^2 ) + x_1 [ (\alpha_6 +\alpha_4)x_2y_2 + (\alpha_3+\alpha_5)x_3y_3 ] &=& \theta y_1 , \\
y_2 ( \alpha_1 x_2^2 + \alpha_4 x_1^2 + \alpha_3 x_3^2 ) + x_2 [ (\alpha_6 +\alpha_4)x_1y_1 + (\alpha_3+\alpha_5)x_3y_3 ] &=& \theta y_2 , \\
y_3 ( \alpha_3 x_1^2 + \alpha_3 x_2^2 + \alpha_2 x_3^2 ) + (\alpha_3+\alpha_5) x_3 ( x_1y_1 + x_2y_2 ) &=& \theta y_3 .
\end{eqnarray*}

Using the constraints $x_1^2+x_2^2+x_3^2 = 1$ and $y_1^2+y_2^2+y_3^2 = 1$, we have the following equivalent reformulation.
\begin{eqnarray*}
x_1 [ ( \alpha_1-\alpha_3) y_1^2 + (\alpha_4-\alpha_3) y_2^2 ] + y_1 [ (\alpha_6 +\alpha_4)x_2y_2 + (\alpha_3+\alpha_5)x_3y_3 ] &=& ( \theta - \alpha_3 ) x_1 , \\
x_2 [ ( \alpha_1-\alpha_3) y_2^2 + (\alpha_4-\alpha_3) y_1^2 ] + y_2 [ (\alpha_6 +\alpha_4)x_1y_1 + (\alpha_3+\alpha_5)x_3y_3 ] &=& ( \theta - \alpha_3 ) x_2 , \\
(\alpha_2-\alpha_3) x_3 y_3^2 + (\alpha_3+\alpha_5) y_3 (x_1y_1 + x_2y_2 ) &=& ( \theta - \alpha_3 ) x_3 , \\
y_1 [ ( \alpha_1-\alpha_3) x_1^2 + ( \alpha_4-\alpha_3) x_2^2 ] + x_1 [ (\alpha_6 +\alpha_4)x_2y_2 + (\alpha_3+\alpha_5)x_3y_3 ] &=& ( \theta - \alpha_3 ) y_1 , \\
y_2 [ ( \alpha_1-\alpha_3) x_2^2 + ( \alpha_4-\alpha_3) x_1^2 ] + x_2 [ (\alpha_6 +\alpha_4)x_1y_1 + (\alpha_3+\alpha_5)x_3y_3 ] &=& ( \theta - \alpha_3 ) y_2 , \\
(\alpha_2-\alpha_3) y_3 x_3^2 + (\alpha_3+\alpha_5) x_3 ( x_1y_1 + x_2y_2 ) &=& ( \theta - \alpha_3 ) y_3 .
\end{eqnarray*}

If we set $\alpha_6+\alpha_4 = \alpha_1-\alpha_4$, we can recover the polynomials arising from the transverse isotropic case.
From the system of polynomial equations above, we can easily verify that $( \theta - \alpha_3 ) x_3^2  = ( \theta - \alpha_3 ) y_3^2$. Hence, we have $\theta = \alpha_3$ or $|x_3|=|y_3|$.

For the case where $x_3 = y_3 =0$, we have
\begin{eqnarray*}
x_1 [ ( \alpha_1-\alpha_3) y_1^2 + (\alpha_4-\alpha_3) y_2^2 ] +  (\alpha_6 +\alpha_4) y_1 x_2y_2   &=& ( \theta - \alpha_3 ) x_1 , \\
x_2 [ ( \alpha_1-\alpha_3) y_2^2 + (\alpha_4-\alpha_3) y_1^2 ] +  (\alpha_6 +\alpha_4) y_2 x_1y_1   &=& ( \theta - \alpha_3 ) x_2 , \\
y_1 [ ( \alpha_1-\alpha_3) x_1^2 + ( \alpha_4-\alpha_3) x_2^2 ] +  (\alpha_6 +\alpha_4) x_1 x_2y_2   &=& ( \theta - \alpha_3 ) y_1 , \\
y_2 [ ( \alpha_1-\alpha_3) x_2^2 + ( \alpha_4-\alpha_3) x_1^2 ] +  (\alpha_6 +\alpha_4) x_2 x_1y_1   &=& ( \theta - \alpha_3 ) y_2   .
\end{eqnarray*}
We obtain $\theta = \alpha_1$, $\alpha_4$, $\frac{\alpha_1-\alpha_6}{2}$ and $\frac{\alpha_1+2\alpha_4+\alpha_6}{2}$.
Then in the following we can only consider $ |x_3| = |y_3| \neq 0$.

(I)
For the case where $x_3 = y_3$, we have
\begin{eqnarray*}
x_1 [ ( \alpha_1-\alpha_3) y_1^2 + (\alpha_4-\alpha_3) y_2^2 ] + y_1 [ (\alpha_6 +\alpha_4)x_2y_2 + (\alpha_3+\alpha_5)x_3^2 ] &=& ( \theta - \alpha_3 ) x_1 , \\
x_2 [ ( \alpha_1-\alpha_3) y_2^2 + (\alpha_4-\alpha_3) y_1^2 ] + y_2 [ (\alpha_6 +\alpha_4)x_1y_1 + (\alpha_3+\alpha_5)x_3^2 ] &=& ( \theta - \alpha_3 ) x_2 , \\
y_1 [ ( \alpha_1-\alpha_3) x_1^2 + ( \alpha_4-\alpha_3) x_2^2 ] + x_1 [ (\alpha_6 +\alpha_4)x_2y_2 + (\alpha_3+\alpha_5)x_3^2 ] &=& ( \theta - \alpha_3 ) y_1 , \\
y_2 [ ( \alpha_1-\alpha_3) x_2^2 + ( \alpha_4-\alpha_3) x_1^2 ] + x_2 [ (\alpha_6 +\alpha_4)x_1y_1 + (\alpha_3+\alpha_5)x_3^2 ] &=& ( \theta - \alpha_3 ) y_2 , \\
(\alpha_2-\alpha_3) x_3^3 + (\alpha_3+\alpha_5) x_3 (x_1y_1 + x_2y_2 ) &=& ( \theta - \alpha_3 ) x_3.
\end{eqnarray*}

%-------------------CANNOT SOLVE EQN!----------------

Using the symmetry of the system of polynomials among $x_1 \leftrightarrow x_2$ and $y_1 \leftrightarrow y_2$, we have
\begin{eqnarray*}
( \alpha_1-2\alpha_3+\alpha_4)x_1  y_1^2  + y_1 [ (\alpha_6 +\alpha_4)x_1y_1 + (\alpha_3+\alpha_5)x_3^2 ] &=& ( \theta - \alpha_3 ) x_1 , \\
( \alpha_1-2\alpha_3+\alpha_4) y_1 x_1^2  + x_1 [ (\alpha_6 +\alpha_4)x_1y_1 + (\alpha_3+\alpha_5)x_3^2 ] &=& ( \theta - \alpha_3 ) y_1 ,\\
(\alpha_2-\alpha_3) x_3^3 + 2(\alpha_3+\alpha_5) x_3 x_1y_1 &=& ( \theta - \alpha_3 ) x_3
\end{eqnarray*}
The we have the M-eigenvalues: $\theta = \alpha_2$, $\theta = \frac{1}{2}(\alpha_6+2\alpha_4+\alpha_1)$,
$\theta = \frac{\alpha_2(\alpha_6+2\alpha_4+\alpha_1) - 2\alpha_5^2}{2\alpha_2 + \alpha_6+2\alpha_4+\alpha_1  +4 \alpha_5}$, $\theta = \frac{(\alpha_6+2\alpha_4+\alpha_1) \alpha_2 - 2(2\alpha_3+\alpha_5)^2}{\alpha_6+2\alpha_4+\alpha_1 + 2\alpha_2 -4(2\alpha_3+\alpha_5)}$.

Using the symmetry: $x_1 \leftrightarrow y_1$, $x_2 \leftrightarrow y_2$, we have
\begin{eqnarray*}
x_1 [ ( \alpha_1-\alpha_3) x_1^2 + (\alpha_4-\alpha_3) x_2^2 ] + x_1 [ (\alpha_6 +\alpha_4)x_2^2 + (\alpha_3+\alpha_5)x_3^2 ] &=& ( \theta - \alpha_3 ) x_1 , \\
x_2 [ ( \alpha_1-\alpha_3) x_2^2 + (\alpha_4-\alpha_3) x_1^2 ] + x_2 [ (\alpha_6 +\alpha_4)x_1^2 + (\alpha_3+\alpha_5)x_3^2 ] &=& ( \theta - \alpha_3 ) x_2 , \\
(\alpha_2-\alpha_3) x_3^3 + (\alpha_3+\alpha_5) x_3 (x_1^2 + x_2^2 ) &=& ( \theta - \alpha_3 ) x_3.
\end{eqnarray*}
And we obtain the M-eigenvalues: $\theta = \alpha_1$, $\theta =\alpha_2$,
$\theta =\frac{1}{2}(\alpha_6+2\alpha_4+\alpha_1)$,
$\theta = \frac{\alpha_1\alpha_2 - (2\alpha_3+\alpha_5)^2}{\alpha_1+\alpha_2 - 2(2\alpha_3+\alpha_5)}$,
$\theta = \frac{(\alpha_6+2\alpha_4+\alpha_1) \alpha_2 - 2(2\alpha_3+\alpha_5)^2}{\alpha_6+2\alpha_4+\alpha_1 + 2\alpha_2 -4(2\alpha_3+\alpha_5)}$.
\footnote{
Finally, again using the symmetry,
we arrive at the problem on three unknowns $x_1, x_3$ and $\theta$.
\begin{eqnarray*}
(\alpha_1-2\alpha_3+2\alpha_4+\alpha_6)  x_1^3 + (\alpha_3+\alpha_5)x_1x_3^2 &=& (\theta - \alpha_3) x_1 , \\
(\alpha_2-2\alpha_3-\alpha_5) x_3^3   &=& (\theta - 2\alpha_3 - \alpha_5) x_3 , \\
2 x_1^2 + x_3^2 &=& 1.
\end{eqnarray*}
Hence, we have the M-eigenvalues $\theta = \alpha_2$,
$\theta =\frac{1}{2}(\alpha_6+2\alpha_4+\alpha_1)$,  and
$\theta = \frac{(\alpha_6+2\alpha_4+\alpha_1) \alpha_2 - 2(2\alpha_3+\alpha_5)^2}{\alpha_6+2\alpha_4+\alpha_1 + 2\alpha_2 -4(2\alpha_3+\alpha_5)}$.
}

(II)
For the case where $x_3 =- y_3$, we have
\begin{eqnarray*}
x_1 [ ( \alpha_1-\alpha_3) y_1^2 + (\alpha_4-\alpha_3) y_2^2 ] + y_1 [ (\alpha_6 +\alpha_4)x_2y_2 - (\alpha_3+\alpha_5)x_3^2 ] &=& ( \theta - \alpha_3 ) x_1 , \\
x_2 [ ( \alpha_1-\alpha_3) y_2^2 + (\alpha_4-\alpha_3) y_1^2 ] + y_2 [ (\alpha_6 +\alpha_4)x_1y_1 - (\alpha_3+\alpha_5)x_3^2 ] &=& ( \theta - \alpha_3 ) x_2 , \\
y_1 [ ( \alpha_1-\alpha_3) x_1^2 + ( \alpha_4-\alpha_3) x_2^2 ] + x_1 [ (\alpha_6 +\alpha_4)x_2y_2 -(\alpha_3+\alpha_5)x_3^2 ] &=& ( \theta - \alpha_3 ) y_1 , \\
y_2 [ ( \alpha_1-\alpha_3) x_2^2 + ( \alpha_4-\alpha_3) x_1^2 ] + x_2 [ (\alpha_6 +\alpha_4)x_1y_1 - (\alpha_3+\alpha_5)x_3^2 ] &=& ( \theta - \alpha_3 ) y_2 , \\
(\alpha_2-\alpha_3) x_3^3 - (\alpha_3+\alpha_5) x_3 (x_1y_1 + x_2y_2 ) &=& ( \theta - \alpha_3 ) x_3.
\end{eqnarray*}

%-------------------CANNOT SOLVE EQN?!----------------

Using the symmetry of the system of polynomials among $x_1 \leftrightarrow x_2$ and $y_1 \leftrightarrow y_2$, we have
\begin{eqnarray*}
( \alpha_1-2\alpha_3+\alpha_4)x_1  y_1^2  + y_1 [ (\alpha_6 +\alpha_4)x_1y_1 - (\alpha_3+\alpha_5)x_3^2 ] &=& ( \theta - \alpha_3 ) x_1 , \\
( \alpha_1-2\alpha_3+\alpha_4) y_1 x_1^2  + x_1 [ (\alpha_6 +\alpha_4)x_1y_1 - (\alpha_3+\alpha_5)x_3^2 ] &=& ( \theta - \alpha_3 ) y_1 ,\\
(\alpha_2-\alpha_3) x_3^3 - 2(\alpha_3+\alpha_5) x_3 x_1y_1 &=& ( \theta - \alpha_3 ) x_3.
\end{eqnarray*}
The corresponding M-eigenvalues are as follows: $\theta = \alpha_2$, $\theta = \frac{1}{2}(\alpha_6+2\alpha_4+\alpha_1)$,  $\theta = \frac{\alpha_2(\alpha_6+2\alpha_4+\alpha_1) - 2\alpha_5^2}{2\alpha_2 + \alpha_6+2\alpha_4+\alpha_1  +4 \alpha_5}$, $\theta = \frac{(\alpha_6+2\alpha_4+\alpha_1) \alpha_2 - 2(2\alpha_3+\alpha_5)^2}{\alpha_6+2\alpha_4+\alpha_1 + 2\alpha_2 -4(2\alpha_3+\alpha_5)}$.

Using the symmetry: $x_1 \leftrightarrow y_1$ and $x_2 \leftrightarrow y_2$, we have
\begin{eqnarray*}
x_1 [ ( \alpha_1-\alpha_3) x_1^2 + (\alpha_4-\alpha_3) x_2^2 ] + x_1 [ (\alpha_6 +\alpha_4)x_2^2 - (\alpha_3+\alpha_5)x_3^2 ] &=& ( \theta - \alpha_3 ) x_1 , \\
x_2 [ ( \alpha_1-\alpha_3) x_2^2 + (\alpha_4-\alpha_3) x_1^2 ] + x_2 [ (\alpha_6 +\alpha_4)x_1^2 - (\alpha_3+\alpha_5)x_3^2 ] &=& ( \theta - \alpha_3 ) x_2 , \\
(\alpha_2-\alpha_3) x_3^3 - (\alpha_3+\alpha_5) x_3 (x_1^2 + x_2^2 ) &=& ( \theta - \alpha_3 ) x_3.
\end{eqnarray*}
And the M-eigenvalues read as follows:
$\theta=\alpha_1$, $\theta=\alpha_2$,
$\theta = \frac{\alpha_1\alpha_2-\alpha_5^2}{\alpha_1+\alpha_2+2\alpha_5}$,
$\theta = \frac{1}{2}(\alpha_6+2\alpha_4+\alpha_1)$,
$\theta = \frac{\alpha_2(\alpha_6+2\alpha_4+\alpha_1) - 2\alpha_5^2}{2\alpha_2 + \alpha_6+2\alpha_4+\alpha_1  +4 \alpha_5}$.
\footnote{
Using more symmetry,
we arrive at the problem on three unknowns $x_1, x_3$ and $\theta$.
\begin{eqnarray*}
(\alpha_1-2\alpha_3+2\alpha_4+\alpha_6)  x_1^3 - (\alpha_3+\alpha_5)x_1x_3^2 &=& (\theta - \alpha_3) x_1 , \\
(\alpha_2+\alpha_5) x_3^3   &=& (\theta + \alpha_5) x_3,\\
2 x_1^2 + x_3^2 &=& 1.
\end{eqnarray*}
Therefore, we have the M-eigenvalues: $\theta = \alpha_2$, $\theta = \frac{1}{2}(\alpha_6+2\alpha_4+\alpha_1)$ and $\theta = \frac{\alpha_2(\alpha_6+2\alpha_4+\alpha_1) - 2\alpha_5^2}{2\alpha_2 + \alpha_6+2\alpha_4+\alpha_1  +4 \alpha_5}$.
}

%----------------------------------------------------

%-*** symmetric polynomials ***-
Hence, solving the polynomial system, we find the following M-eigenvalues:
%$$\theta_1 = \alpha_1, \qquad \theta_2 = \alpha_2, \qquad \theta_3 = \alpha_3, \qquad \theta_4=\alpha_4, $$
\begin{eqnarray*}
%\theta_1 &=& \alpha_1, \\
%\theta_2 &=& \alpha_2, \\
%\theta_3 &=& \alpha_3, \\
%\theta_4 &=& \alpha_4, \\
\theta_i &=& \alpha_i, \quad (i = 1,2,3,4) \\
\theta_5 &=& \frac{\alpha_1\alpha_2-\alpha_5^2}{\alpha_1+\alpha_2+2\alpha_5}, \\
\theta_6 &=& \frac{\alpha_1\alpha_2 - (2\alpha_3+\alpha_5)^2}{\alpha_1+\alpha_2 - 2(2\alpha_3+\alpha_5)},\\
\theta_7 &=& \frac{1}{2}(\alpha_6+2\alpha_4+\alpha_1), \\
\theta_8 &=& \frac{1}{2} (\alpha_1-\alpha_6), \\
\theta_9 &=& \frac{\alpha_2(\alpha_6+2\alpha_4+\alpha_1) - 2\alpha_5^2}{2\alpha_2 + \alpha_6+2\alpha_4+\alpha_1  +4 \alpha_5}, \\
\theta_{10} &=& \frac{(\alpha_6+2\alpha_4+\alpha_1) \alpha_2 - 2(2\alpha_3+\alpha_5)^2}{\alpha_6+2\alpha_4+\alpha_1 + 2\alpha_2 -4(2\alpha_3+\alpha_5)}.
\end{eqnarray*}
Note that  $\theta_9 = \frac{ \theta_7 \alpha_2 - \alpha_5^2}{\theta7 + \alpha_2 +2\alpha_5}$, $\theta_{10} = \frac{ \theta_7 \alpha_2 - (2\alpha_3+\alpha_5)^2}{\theta7 + \alpha_2 -2(2\alpha_3+\alpha_5)}$.
When $\alpha_6=\alpha_1 - 2\alpha_4$, we have $\theta_7 = \theta_1$, $\theta_8 = \theta_4$, $\theta_9=\theta_5$ and $\theta_{10} = \theta_6$, and we recover the M-eigenvalues for the polar anisotropic case.

The corresponding Z-eigenvalue problem reads
\begin{eqnarray*}
x_1 [ ( \alpha_1-2\alpha_3-\alpha_5) x_1^2 + (\alpha_6+2\alpha_4-2\alpha_3-\alpha_5) x_2^2 ]
&=& ( \eta - 2\alpha_3 -\alpha_5 ) x_1 , \\
x_2 [ ( \alpha_1-2\alpha_3-\alpha_5) x_2^2 + (\alpha_6+2\alpha_4-2\alpha_3-\alpha_5) x_1^2 ]  &=& ( \eta - 2\alpha_3 -\alpha_5 ) x_2 , \\
( \alpha_2 - 2\alpha_3 -\alpha_5 ) x_3^3  &=& ( \eta - 2\alpha_3 -\alpha_5 ) x_3  .
\end{eqnarray*}
We can calculate the Z-eigenvalues:
\begin{eqnarray*}
\eta_1 &=& \frac{1}{2}(\alpha_6+2\alpha_4+\alpha_1), \\
\eta_2 &=& \alpha_2, \\
\eta_3 &=& \frac{\alpha_1\alpha_2 - (2\alpha_3+\alpha_5)^2}{\alpha_1+\alpha_2 - 2(2\alpha_3+\alpha_5)},\\
\eta_4 &=& \frac{(\alpha_6+2\alpha_4+\alpha_1) \alpha_2 - 2(2\alpha_3+\alpha_5)^2}{\alpha_6+2\alpha_4+\alpha_1 + 2\alpha_2 -4(2\alpha_3+\alpha_5)}.
\end{eqnarray*}

If $\alpha_6=\alpha_1 - 2\alpha_4$, then $\eta_1 = \alpha_1$ and $\eta_3 = \eta_4$, and the results reduce to the Z-eigenvalues of the polar anisotropic system.

%Define $\zeta:= \frac{1}{2}(\alpha_6+2\alpha_4+\alpha_1)$.
We can easily check that the necessary and sufficient conditions for $\theta_i >0 (i=1,\cdots, 10)$ are the following inequalities.
\begin{equation}\label{Tetragonal1}
\alpha_i > 0 ~ (i=1,2,3,4), \qquad
\alpha_6 + 2 \alpha_4 + \alpha_1 > 0, \qquad
\alpha_1 > \alpha_6 ,
\end{equation}
\begin{equation}\label{Tetragonal2}
| \alpha_5 | < \sqrt{\alpha_1 \alpha_2}  , \qquad
| \alpha_5 | < \sqrt{ \theta_7 \alpha_2 } ,
\end{equation}
\begin{equation}\label{Tetragonal3A}
| 2 \alpha_3 + \alpha_5 | < \sqrt{\alpha_1 \alpha_2}  \quad
\text{or} \quad
 2 \alpha_3 + \alpha_5 > \frac{1}{2}(\alpha_1+\alpha_2) ,
\end{equation}

\begin{equation}\label{Tetragonal3B}
| 2 \alpha_3 + \alpha_5 | < \sqrt{ \theta_7 \alpha_2 }  \quad
\text{or} \quad 2 \alpha_3 + \alpha_5 > \frac{1}{2}(\theta_7 + \alpha_2) .
\end{equation}

Using the definitions of $\alpha_i~(i=1,\cdots,6)$, the six inequalities in \eqref{Tetragonal1} and the two inequalities in \eqref{Tetragonal2} can be equivalently expressed as follows.
\begin{equation}\label{eqn:Chirita_tetragonal0}
 c_{11} >0, \quad c_{33}>0, \quad c_{55} >0, \quad c_{66}>0, \quad c_{12}+2c_{66}+c_{11}>0, \quad c_{11}>c_{12} ,
\end{equation}
\begin{equation}\label{eqn:Tetragonal2inConstant}
| c_{13} | < \min \{ \sqrt{c_{11}c_{33}}, \sqrt{\zeta c_{33}} \}  ,
\end{equation}
where  $\zeta=  \frac{1}{2}(c_{12}+2c_{66}+c_{11})  $.
The inequalities in \eqref{Tetragonal3A} and \eqref{Tetragonal3B} can be reexpressed as follows.
\begin{equation}\label{eqn:TetragonalSufficient1}
|2c_{55}+c_{13}| < \min \{ \sqrt{c_{11}c_{33}}, \sqrt{\zeta c_{33}} \},
\end{equation}
\begin{equation}\label{eqn:TetragonalSufficient2}
2c_{55}+c_{13} > \frac{1}{2} \left( c_{33} + \max\{c_{11}, \zeta\} \right) ,
\end{equation}
\begin{equation}\label{eqn:TetragonalSufficient3}
\frac{1}{2} \left( \zeta + c_{33}  \right) < 2c_{55}+c_{13} <   \sqrt{c_{11}c_{33}},
\end{equation}
or
\begin{equation}\label{eqn:TetragonalSufficient4}
\frac{1}{2} \left( c_{11} + c_{33} \right) < 2c_{55}+c_{13} <  \sqrt{\zeta c_{33}}.
\end{equation}
In fact, we find that the elastic constants of the materials TeO$_2$, Ba$_2$Si$_2$TiO$_8$,
Ca$_{10}$Mg$_2$Al$_4$(SiO$_4$)$_5$(Si$_2$O$_7$)$_2$(OH)$_4$, Paratellurite, Scapolite, Fresnoite, etc. satisfy \eqref{eqn:TetragonalSufficient1}; the elastic constants of SiO$_2$ and Zircon satisfy \eqref{eqn:TetragonalSufficient2}; the elastic constants of TiO$_2$, GeO$_2$ and SnO$_2$ satisfy \eqref{eqn:TetragonalSufficient4}; the elastic constants of ZrSiO$_4$ satisfy \eqref{eqn:TetragonalSufficient3} \cite[P.50]{Ahrens_Book1995}.

%That is,
%\begin{equation}
%\max \{ | c_{33} |, |2c_{55}+c_{13}| \} < \min \{ \sqrt{c_{11}c_{33}}, \sqrt{\zeta} \}.
%\end{equation}

{\bf Remarks.}
From the inequalities in \eqref{eqn:Tetragonal2inConstant} and \eqref{eqn:TetragonalSufficient1}, we have
\begin{equation}\label{eqn:Chirita_tetragonal1}
 -2 c_{55} - \sqrt{c_{11}c_{33}} < c_{13} < \sqrt{c_{11}c_{33}} .
\end{equation}
Since $| 2 c_{55} + c_{13} | \geq | | c_{55} + c_{13} | - c_{55}  |$, from the inequalities in \eqref{eqn:TetragonalSufficient1} we obtain
\begin{equation}\label{eqn:Chirita_tetragonal2}
( | c_{55} + c_{13} | - c_{55})^2 < \frac{1}{2} c_{33} \left( c_{11} + c_{12} + \min\{ c_{11}-c_{12}, 2 c_{66} \} \right) .
\end{equation}
The strong ellipticity conditions \eqref{eqn:Chirita_tetragonal0}, \eqref{eqn:Chirita_tetragonal1} and \eqref{eqn:Chirita_tetragonal2} already appeared in \cite{SChirita_IJSS08} (see \cite[(1.6)-(1.8) of Theorem 1]{SChirita_IJSS08}).

Using the fact that $| 2 c_{55} + c_{13} | \geq | c_{55} + c_{13} | - c_{55}$ together with the definition $\zeta$, and the inequalities in \eqref{eqn:TetragonalSufficient1},  we can derive the following inequalities.
$$|c_{55}+c_{13}| < c_{55} + \sqrt{c_{11}c_{33}},  \qquad |c_{55}+c_{13}| < c_{55} + \sqrt{\zeta c_{33}}. $$
That is, $|c_{55}+c_{13}| < c_{55} + \min \{ \sqrt{c_{11}c_{33}}, \sqrt{\zeta c_{33}} \} $. Thus we recover the conditions given in \cite{ZubovRudev16} (see \cite[(2.4)-(2.6)]{ZubovRudev16}).
%From the last three inequalities, we can easily derive that ... includes (3.4) in ...

\section{Orthotropy}

An orthotropic material has three mutually perpendicular symmetry planes.   With basis vectors perpendicular to the symmetry planes, the elastic stiffness matrix has the  following form:
$$ ( c_{\alpha \beta } ) = {\begin{bmatrix}c_{11} &c_{12} & c_{13}&0&0&0 \\c_{12}&c_{22}&c_{23}&0&0&0 \\ c_{13}&c_{23}&c_{33}&0&0&0 \\0&0&0&c_{44}&0&0\\0&0&0&0&c_{55}&0\\0&0&0&0&0&c_{66}\end{bmatrix}}.$$
The elasticity tensor is described by nine independent material parameters:
$c_{11}= C_{1111} :=\delta_1$, $c_{22} = C_{2222} :=\delta_2$, $c_{33}=C_{3333} := \delta_3$, $c_{44}=C_{2323}=C_{3232} :=\delta_4$, $c_{55}=C_{1313}=C_{3131} :=\delta_5$, $c_{66}=C_{1212}=C_{2121} :=\delta_6$, $c_{23}=C_{2233} := \delta_7$, $c_{13}=C_{1133} := \delta_8$, $c_{12}=C_{1122} :=\delta_9$. The components of the elastic stiffness matrix are related to the elastic constants \cite[p.83]{Bower_book}\cite[p.47]{PaoloVannucci_book18}.

The corresponding M-eigenvalue problem reads as follows.
\begin{eqnarray*}
x_1 ( \delta_1 y_1^2 + \delta_6 y_2^2 + \delta_5 y_3^2 ) + y_1 [ (\delta_9 +\delta_6)x_2y_2 + (\delta_8+\delta_5)x_3y_3 ] &=& \theta x_1 , \\
x_2 ( \delta_6 y_1^2 + \delta_2 y_2^2 + \delta_4 y_3^2 ) + y_2 [ (\delta_7+\delta_4)x_3y_3 + (\delta_9 +\delta_6)x_1y_1 ] &=& \theta x_2 , \\
x_3 ( \delta_5 y_1^2 + \delta_4 y_2^2 + \delta_3 y_3^2 ) + y_3 [ (\delta_8+\delta_5) x_1y_1 + (\delta_7+\delta_4) x_2y_2 ] &=& \theta x_3 , \\
y_1 ( \delta_1 x_1^2 + \delta_6 x_2^2 + \delta_5 x_3^2 ) + x_1 [ (\delta_9 +\delta_6)x_2y_2 + (\delta_8+\delta_5)x_3y_3 ] &=& \theta y_1 , \\
y_2 ( \delta_6 x_1^2 + \delta_2 x_2^2 + \delta_4 x_3^2 ) + x_2 [ (\delta_7+\delta_4)x_3y_3 + (\delta_9 +\delta_6)x_1y_1 ] &=& \theta y_2 , \\
y_3 ( \delta_5 x_1^2 + \delta_4 x_2^2 + \delta_3 x_3^2 ) + x_3 [ (\delta_8+\delta_5) x_1y_1 + (\delta_7+\delta_4) x_2y_2 ] &=& \theta y_3 .
\end{eqnarray*}

%-------------------CANNOT SOLVE EQN !!!----------------

\footnote{ In principal we can solve it via resultants \cite{Book2005CLO}. If we use the symmetry: $x_1 \leftrightarrow \pm y_1$, $x_2 \leftrightarrow \pm y_2$, $x_3 \leftrightarrow \pm y_3$. Setting $x= \pm y$, we have
\begin{eqnarray*}
x_1 [ \delta_1 x_1^2  +  (\delta_9 +2\delta_6)x_2^2 + (\delta_8+2\delta_5)x_3^2 ] &=& \theta x_1 , \\
x_2 [ \delta_2 x_2^2  +  (\delta_7+2\delta_4)x_3^2  + (\delta_9 +2\delta_6)x_1^2 ] &=& \theta x_2 , \\
x_3 [ \delta_3 x_3^2  +  (\delta_8 +2\delta_5) x_1^2 + (\delta_7+2\delta_4) x_2^2 ] &=& \theta x_3  .
\end{eqnarray*}
It yields $\theta_1, \theta_2, \theta_3, \theta_7, \theta_8, \theta_9, \theta_{13}$. The special cases where $x_j=y_j=0$ yield other $\theta$'s.
}

We can find the M-eigenvalues:
\begin{eqnarray*}
&&\theta_i = \delta_i, \quad (i = 1,\cdots,6) \\
&&\theta_7 = \frac{\delta_2\delta_3-\delta_7^2}{\delta_2+\delta_3+2\delta_7}, \\
&&\theta_8 = \frac{\delta_1\delta_3-\delta_8^2}{\delta_1+\delta_3+2\delta_8}, \\
&&\theta_9 = \frac{\delta_1\delta_2-\delta_9^2}{\delta_1+\delta_2+2\delta_9}, \\
&&\theta_{10} = \frac{\delta_2\delta_3 - (2\delta_4+\delta_7)^2}{\delta_2+\delta_3 - 2(2\delta_4+\delta_7)},\\
&&\theta_{11} = \frac{\delta_1\delta_3 - (2\delta_5+\delta_8)^2}{\delta_1+\delta_3 - 2(2\delta_5+\delta_8)},\\
&&\theta_{12} = \frac{\delta_1\delta_2 - (2\delta_6+\delta_9)^2}{\delta_1+\delta_2 - 2(2\delta_6+\delta_9)},
\end{eqnarray*}
and
\begin{eqnarray*}
&& \theta_{13} = ( 4 \delta_1 \delta_4^2 + 4 \delta_2 \delta_5^2  + 4 \delta_3 \delta_6^2 + \delta_1 \delta_7^2 + \delta_2 \delta_8^2 + \delta_3 \delta_9^2  -\delta_1 \delta_2 \delta_3 - 16 \delta_4 \delta_5 \delta_6  \\
&& - 8 \delta_5 \delta_6 \delta_7 - 4 \delta_6 \delta_7 \delta_8 - 2 \delta_7 \delta_8 \delta_9 + 4 \delta_1 \delta_4 \delta_7 + 4 \delta_2 \delta_5 \delta_8  + 4 \delta_3 \delta_6 \delta_9 - 8 \delta_4 \delta_6 \delta_8 \\
&& - 8 \delta_4 \delta_5 \delta_9 - 4 \delta_5 \delta_7 \delta_9 - 4 \delta_4 \delta_8 \delta_9  ) %
/ ( 4 \delta_4^2+ 4 \delta_5^2+ 4 \delta_6^2+ \delta_7^2 + \delta_8^2 + \delta_9^2 -\delta_1 \delta_2 \\
&& - \delta_1 \delta_3 - \delta_2 \delta_3 - 8 \delta_4 \delta_5 - 8 \delta_4 \delta_6-   8 \delta_5 \delta_6 - 2 \delta_7 \delta_8 - 2 \delta_7 \delta_9 - 2 \delta_8 \delta_9 + 4 \delta_1 \delta_4 \\
&& + 4 \delta_2 \delta_5 + 4 \delta_3 \delta_6 + 2 \delta_1 \delta_7 + 2 \delta_2 \delta_8  + 2 \delta_3 \delta_9 + 4 \delta_4 \delta_7 - 4 \delta_4 \delta_8 - 4 \delta_4 \delta_9
- 4 \delta_5 \delta_7 \\
&& + 4 \delta_5 \delta_8 - 4 \delta_5 \delta_9 - 4 \delta_6 \delta_7 - 4 \delta_6 \delta_8 + 4 \delta_6 \delta_9  ).
\end{eqnarray*}
For the M-eigenvalue $\theta_{13}$, we can verify that
\begin{eqnarray*}
\text{the numerator} =  \frac{1}{2}[(\delta_9 + 2\delta_6 + \delta_2)\delta_3 - 2(\delta_7+2\delta_4)^2] (\delta_9 + 2\delta_6 - \delta_1) \\
+ \frac{1}{2}[(\delta_9 + 2\delta_6 + \delta_1)\delta_3 - 2(\delta_8+2\delta_5)^2] (\delta_9 + 2\delta_6 - \delta_2) \\
+ (\delta_7-\delta_8+2\delta_4-2\delta_5)^2 (\delta_9 + 2\delta_6),
\end{eqnarray*}
\begin{eqnarray*}
\text{the denominator} = \frac{1}{2}[(\delta_9 + 2\delta_6 + \delta_2 + 2\delta_3 - 4(\delta_7+2\delta_4)] (\delta_9 + 2\delta_6 - \delta_1) \\
+ \frac{1}{2}[(\delta_9 + 2\delta_6 + \delta_1 + 2\delta_3 - 4(\delta_8+2\delta_5)] (\delta_9 + 2\delta_6 - \delta_2) \\
+ (\delta_7-\delta_8+2\delta_4-2\delta_5)^2  .
\end{eqnarray*}
These expressions will help us to reduce $\theta_{13}$ to the M-eigenvalue of tetragonal case and provide a sufficient condition for $\theta_{13}$ being positive.

When $\delta_1=\delta_2=\alpha_1$, $\delta_3=\alpha_2$, $\delta_4=\delta_5=\alpha_3$,
$\delta_6=\alpha_4$, $\delta_7=\delta_8=\alpha_5$ and $\delta_9=\alpha_6$, we have
$$\theta^{(O)}_1=\theta^{(O)}_2 = \theta^{(T)}_1, \quad
\theta^{(O)}_3 = \theta^{(T)}_2, \quad
\theta^{(O)}_4=\theta^{(O)}_5 = \theta^{(T)}_3, \quad
\theta^{(O)}_6 = \theta^{(T)}_4, $$
$$
\theta^{(O)}_7 = \theta^{(O)}_8 = \theta^{(T)}_5, \quad
\theta^{(O)}_9 = \theta^{(T)}_8 , \quad
\theta^{(O)}_{10}=\theta^{(O)}_{11} = \theta^{(T)}_6,\quad
\theta^{(O)}_{12} = \theta^{(T)}_7, $$
and
%\begin{eqnarray*}
%&& \theta^{(O)}_{13} = \\
%&& (\alpha_1 \alpha_1 \alpha_2 + 16 \alpha_3 \alpha_3 \alpha_4 + 16 \alpha_3 \alpha_4 \alpha_5 + 4 \alpha_4 \alpha_5 \alpha_5 + 2 \alpha_5 \alpha_5 \alpha_6 - 8 \alpha_1 \alpha_3 \alpha_5\\
%&& - 4 \alpha_2 \alpha_4 \alpha_6  + 8 \alpha_3 \alpha_3 \alpha_6 + 8 \alpha_3 \alpha_5 \alpha_6 - 8 \alpha_1 \alpha_3^2  - 4 \alpha_2 \alpha_4^2 - 2 \alpha_1 \alpha_5^2 - \alpha_2 \alpha_6^2 ) \\
%&& / ( 2 \alpha_1 \alpha_2   + 16 \alpha_3 \alpha_4  + 4 \alpha_5 \alpha_6 - 8 \alpha_1 \alpha_3  - 4 \alpha_2 \alpha_4  - 4 \alpha_1 \alpha_5 \\
%&& - 2 \alpha_2 \alpha_6  + 8 \alpha_3 \alpha_6 + 8 \alpha_4 \alpha_5  - 4 \alpha_4 \alpha_6 +\alpha_1^2  - 4 \alpha_4^2 - \alpha_6^2 ).
%\end{eqnarray*}
$$\theta^{(O)}_{13} = \frac{[(\alpha_6+2\alpha_4+\alpha_1) \alpha_2 - 2(2\alpha_3+\alpha_5)^2] (\alpha_6+2\alpha_4-\alpha_1)}{[\alpha_6+2\alpha_4+\alpha_1 + 2\alpha_2 -4(2\alpha_3+\alpha_5)] (\alpha_6+2\alpha_4-\alpha_1)} =\theta^{(T)}_{10} ,$$
where the superscripts $(O)$ and $(T)$ stand for the orthotropic and tetragonal systems respectively. Thus, we recover the M-eigenvalues of the tetragonal system.

Similar to the polar anisotropic case and the tetragonal case, we derive the conditions that allow $\theta_i > 0 ~(i=1,\cdots, 12)$:
\begin{eqnarray*}
\delta_i > 0 ~ (i=1,\cdots,6),  \\
\sqrt{\delta_1 \delta_2} > | \delta_9 |, ~\quad \sqrt{\delta_2 \delta_3} > |\delta_7|, ~\quad \sqrt{\delta_1 \delta_3} > |\delta_8|, \\
\sqrt{\delta_2 \delta_3} > | \delta_7 + 2\delta_4 | \quad \text{or} \quad \delta_7 + 2\delta_4 > \frac{1}{2}( \delta_2 + \delta_3 ), \\
\sqrt{\delta_1 \delta_3} > | \delta_8 + 2\delta_5 | \quad \text{or} \quad \delta_8 + 2\delta_5 > \frac{1}{2}( \delta_1 + \delta_3 )  , \\
\sqrt{\delta_1 \delta_2} > | \delta_9 + 2\delta_6 | \quad \text{or} \quad \delta_9 + 2\delta_6 > \frac{1}{2}( \delta_1 + \delta_2 )  .
\end{eqnarray*}
Suppose that $\sqrt{\delta_2 \delta_3} > | \delta_7 + 2\delta_4 |$ and $\sqrt{\delta_1 \delta_3} > | \delta_8 + 2\delta_5 |$, which guarantee that $\theta_{10}>0$ and $\theta_{11}>0$ respectively.  If $\delta_9 + 2\delta_6 > \max \{ \delta_1, \delta_2 \} >0$, then we have
$\delta_9 + 2\delta_6 - \delta_1>0$,  $\delta_9 + 2\delta_6 - \delta_2>0$,
$\sqrt{ \frac{1}{2}(\delta_9 + 2\delta_6 + \delta_2) \delta_3} >  \sqrt{\delta_2 \delta_3} > |\delta_7 + 2\delta_4|$,
and
$\sqrt{ \frac{1}{2}(\delta_9 + 2\delta_6 + \delta_1) \delta_3} >  \sqrt{\delta_1 \delta_3} > |\delta_8 + 2\delta_5|$,
and hence, $\theta_{13}>0$. Besides, it is naturally to show that $\theta_{12} > 0$ since $2\delta_6 + \delta_9>(\delta_1 + \delta_2)/2$ and $2\delta_6 + \delta_9 > \sqrt{\delta_1 \delta_2}$.

Using the definitions of $\delta_i ~(i=1,\cdots, 9)$ by elasticity constants, we reformulate the above analysis and provide the sufficient conditions for strong ellipticity as follows.
\begin{eqnarray*}
&&c_{ii} > 0 ~ (i=1,\cdots,6),  \\
&&\sqrt{c_{11} c_{22} } > | c_{12} |, \\
&&c_{12} + 2 c_{66} > \max \{ c_{11}, c_{22} \},  \\
&&\sqrt{c_{22} c_{33} } > \max \{ |c_{23}|,  | c_{23} + 2 c_{44} | \} , \\
&&\sqrt{c_{11} c_{33} } > \max \{ |c_{13}|,  | c_{13} + 2 c_{55} | \}.
\end{eqnarray*}
Note that $c_{ii} > 0 ~ (i=1,4,5,6)$ and $\sqrt{c_{11} c_{22} } > | c_{12} |$ are also the necessary conditions for the positive definiteness of elasticity tensor \cite[P.50]{PaoloVannucci_book18}.

The corresponding Z-eigenvalue problem reads
\begin{eqnarray*}
x_1 [ \delta_1 x_1^2 + (\delta_9 +2\delta_6)x_2^2 + (\delta_8+2\delta_5)x_3^2 ] &=& \eta x_1 , \\
x_2 [  \delta_2 x_2^2 + (\delta_9 +2\delta_6)x_1^2 + (\delta_7+2\delta_4)x_3^2 ] &=& \eta x_2 , \\
x_3 [ \delta_3 x_3^2 + (\delta_8+2\delta_5) x_1^2 + (\delta_7+2\delta_4)x_2^2 ] &=& \eta x_3 .
\end{eqnarray*}
We can calculate seven Z-eigenvalues: $\eta = \theta_1$, $\theta_2$, $\theta_3$, $\theta_{10}$, $\theta_{11}$, $\theta_{12}$, and $\theta_{13}$. And the results can reduce to the Z-eigenvalues of the tetragonal system.

\section{Concluding Remarks}

In this note, we show that M-eigenvalues coincide, for example, the bulk modulus $K$ and the shear modulus $G$ in the isotropic case. In fact, all the elastic constants in the diagonal of stiffness matrix are M-eigenvalues, and other M-eigenvalues are the simple linear combinations or the rational fractions of the elastic constants. Using the positiveness condition of the M-eigenvalues, we provide the conditions for strong ellipticity, and recover the existing conditions for strong ellipticity in the cubic case, the polar anisotropic case and the tetragonal case.  We further develop the M-eigenvalue technique to the orthotropic case, and give new sufficient conditions for strong ellipticity of orthotropic case.  This further reveals the physical meanings of the M-eigenvalues and shows that positiveness condition of the M-eigenvalues are more general for strong ellipticity.

\end{document}